\documentclass[12pt]{iopart}

\usepackage{iopams}

\usepackage{graphicx}
\usepackage{epstopdf}

\expandafter\let\csname equation*\endcsname\relax
\expandafter\let\csname endequation*\endcsname\relax
\usepackage{amsmath}
\usepackage{amsthm}
\usepackage{subcaption}
\usepackage[american]{babel}
\usepackage[latin1]{inputenc}
\usepackage{txfonts}
\usepackage{xcolor}
\usepackage[all]{xy}

\usepackage{todonotes}

\usepackage{comment}

\usepackage{hyperref}

\newtheorem{theorem}{Theorem}[section]
\newtheorem{proposition}[theorem]{Proposition}
\newtheorem{lemma}      [theorem]{Lemma}

\newtheorem{definition} [theorem]{Definition}
\newtheorem{hypothesis} [theorem]{Hypothesis}
\newtheorem{remark}[theorem]{Remark}

\usepackage{mathtools}
\DeclarePairedDelimiter\abs {\lvert}{\rvert}
\DeclarePairedDelimiter\norm{\lVert}{\rVert}
\DeclarePairedDelimiterX\inner[2]{\langle}{\rangle}{#1 , #2}

\newcommand{\N}{\mathbb{N}} 
\newcommand{\R}{\mathbb{R}} 
\newcommand{\C}{\mathbb{C}} 
\newcommand{\D}{\mathbb{D}} 
\newcommand{\T}{\mathbb{T}} 

\begin{document}

\title[Chimera states]{Chimera states through invariant manifold theory}

\author{Jaap Eldering$^1$, Jeroen S.W. Lamb$^2$, Tiago Pereira$^{1,2}$, and\\ 
Edmilson Roque dos Santos$^1$
}
\address{$^1$Institute of Mathematical and Computer Sciences, University of S\~ao Paulo, Brazil \\
$^2$Department of Mathematics, Imperial College London, London SW7 2AZ, UK \\
}
\ead{tiago.pereira@imperial.ac.uk}
\begin{abstract}
We establish the existence of chimera states, simultaneously supporting synchronous and asynchronous dynamics, in a network consisting of two symmetrically linked star subnetworks consisting of identical oscillators with shear and Kuramoto--Sakaguchi coupling. We show that the chimera states may be metastable or asymptotically stable. If the intra-star coupling strength is of order $\varepsilon$, the chimera states persist on time scales at least of order $1/\varepsilon$ in general, and on time-scales at least of order $1/\varepsilon^2$ if the intra-star coupling is of Kuramoto--Sakaguchi type. If the intra-star coupling configuration is sparse, the chimeras are asymptotically stable. The analysis relies on a combination of dimensional reduction using a M\"obius symmetry group and techniques from averaging theory and normal hyperbolicity.
\end{abstract}

\maketitle

\tableofcontents

\section{Motivation and main results}
\label{sec:motivation_and_main_results}
In 2002, Kuramoto and Battogtokh~\cite{kuramoto2002} observed the coexistence of spatiotemporal synchronous and asynchronous oscillations in a ring of identical coupled oscillators. This phenomenon of partial synchronisation in networks of identical coupled oscillators, was subsequently branded \textit{chimera} and henceforth observed in a wide range of experimental settings, including lasers~\cite{Larger2015}, photoelectrochemical oscillators~\cite{Haugland2015, Totz_2018_spiral_wave_chimeras} and coupled metronomes~\cite{martens2013}. For comprehensive reviews of the prolific occurrence of  chimeras in physical and numerical experiments, see \cite{panaggio2015,Kementh_a_classification_2016}.

On the theoretical side, various mechanisms for the occurrence of chimeras have been  proposed~\cite{panaggio2015,omel2013coherence,wolfrum2011spectral,abrams2008}. In the limit of infinitely many oscillators, a spectral theory for chimera  states was developed \cite{Omel_chenko_2018}, Notably,  the comprehension of chimera states in large but finite size networks is more challenging~\cite{MATIAS}. Carefully designed coupling functions have been shown to generate chimeras in finite size systems~\cite{Ashwin2015,Bick2016}. Overall, despite a good two decades of fascination with chimeras, we remain far from a comprehensive mathematical understanding of their occurrence and behaviour.

The main objective of this paper is to obtain mathematical results on chimera behaviour in a large but finite size network of identical oscillators, with a specific modular structure and a relatively generic type of coupling. To our best knowledge, our results are the first example of rigorous results in this context. The example concerns the union of two identical, symmetrically coupled star subnetworks.  The key ingredients that underlie our analysis are a dimensional reduction due to a M\"obius group symmetry \cite{watanabe1994},  averaging theory \cite{Sanders2007Averaging} and normally hyperbolic invariant manifold theory \cite{fenichel1971,HPS1977}.

\begin{figure}
         \centering
         \includegraphics[width=0.8\textwidth]{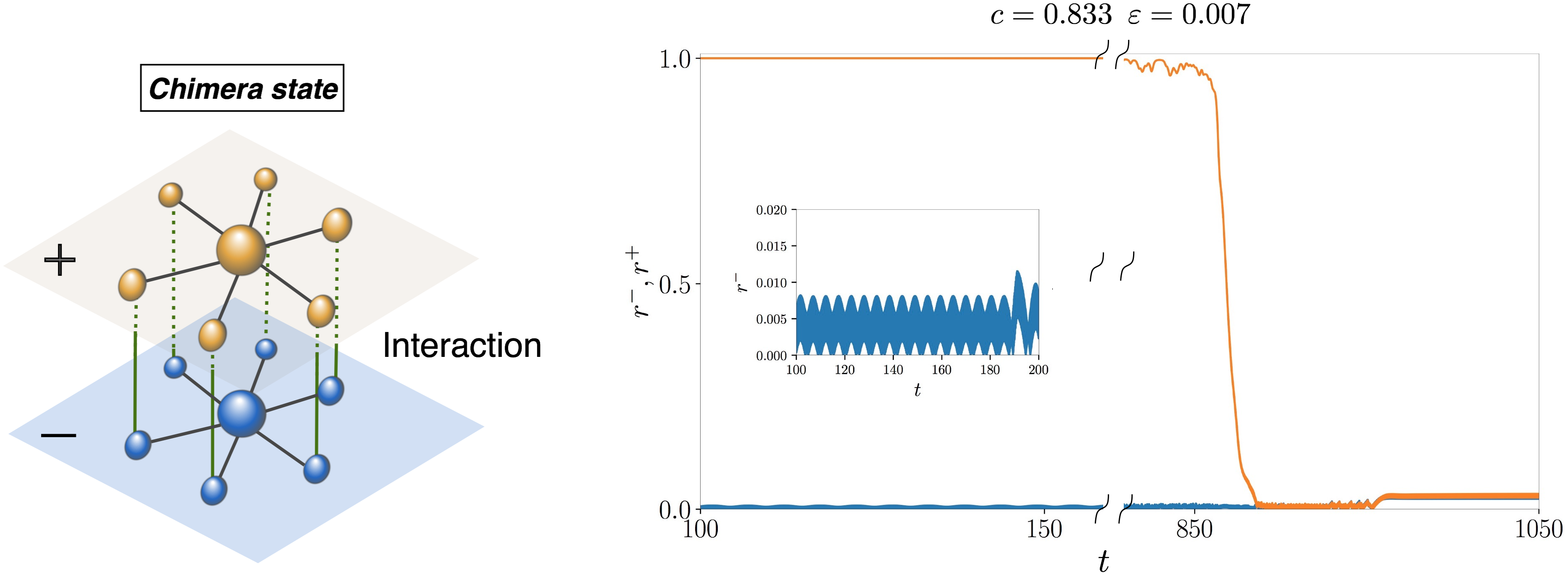}    \caption{Metastability of the chimera states in symmetrically coupled star networks. In the left panel we show a sketch of the network as a union of two symmetrically coupled star subnetworks. In the right panel we show the time series of order parameters $r^{+}$ and $r^{-}$ for two the  star subnetworks $+$ and $-$. Equations of motions are given by~\eqref{eq:coupled_stars_initial} with $N = 200$, coupling between stars as $h(\varphi_j^{\pm}, \varphi_i^{\mp}) = \sin(\varphi_j^{\pm} - \varphi_i^{\mp} + \delta)$ for $i, j \in \{0, \dots N\}$, $\delta = 0.3$, $\omega = 1$ and $ \lambda = 1$, where $c$ is the shear parameter and $\varepsilon$ the intra-star coupling strength.  After approximately $10^3$ cycles of natural frequency the chimera collapses to asynchrounous state. Choosing different parameter values, the chimera may collapse to a synchronous or asynchronous state. For more details, see Section~\ref{numerics}. 
}
        \label{fig:coupled_stars}
\end{figure}

The network we consider is sketched in Figure~\ref{fig:coupled_stars}. It consists of two symmetrically coupled star subnetworks. Each such star subnetwork consists of a central \textit{hub} connected to a large number ($N$) of \textit{leaves}, where all nodes are occupied by identical phase oscillators and interactions arise through a diffusive Kuramoto--Sakaguchi coupling with shear.  The strength of the coupling between nodes in each of the star motif subnetworks is uniform, while the strength of the interaction between the star motif components is represented by another parameter $\varepsilon$. For more details on the equations of motion, see Section~\ref{subsec:model} below.

An important feature of the system is that when isolated, i.e., when $\varepsilon=0$, each of the star motif components exhibits a bistable behaviour of its collective dynamics, exhibiting the coexistence of one stable state in which has all leaves are in \textit{synchrony} and another stable state in which the motion is \textit{asynchronous}. In our setting the bi-stability is induced by the shear. We express the level of synchrony in terms of a complex order parameter $z$, see Eq.~\eqref{eq:order_parameter}, related to the average phase difference between corresponding phase oscillators in the two-star motif subnetworks, where $r=1$ refers to full synchrony and $r=0$ represents full asynchrony.

We consider the dynamics of the coupled subnetworks ($\varepsilon>0$) with shear parameter chosen in the bi-stability regime of the isolated dynamics and initial conditions for each subnetwork are chosen near the alternative stable states, so that the resulting dynamics in the absence of coupling ($\varepsilon=0$) would display a combination of synchronous and asynchronous dynamics. In the presence of coupling ($\varepsilon>0$), numerical experiments show that such states tend to persist (to the limit of our simulation time, approximately $10^6$ cycles of natural frequency) when the coupling strength $\varepsilon$ between the star motif is small, or collapse to a completely synchronized or asynchronous state when this coupling strength is not so small. Representative examples from our simulations are presented in Figure~\ref{fig:coupled_stars}.

In this paper we address mathematically the explanation of these observations for sufficiently large $N$ and sufficiently small $\varepsilon$.  We obtain the following main results:
\begin{description}
\item[{\it (Short metastability)}] Chimera states exist for a time $\mathcal{O}(\varepsilon^{-1})$ \footnote{Here $\mathcal{O}(\cdot)$ stands for Landau's symbols for order functions. So, $t = \mathcal{O}\bigl(\delta(\varepsilon)\bigr)$ is: there exists $\varepsilon_0 \geq 0$ and $K \geq 0$ such that $0 \leq t \leq K |\delta(\varepsilon)|$ for $0 \leq \varepsilon \leq \varepsilon_0$.}, see
Theorem~\ref{thm:chimera-state}~[\ref{enum:thm-eps-time}].
\item[{\it (Long metastability)}] If the intra-star coupling is of Kuramoto--Sakaguchi type
then chimera states exist up to a time $\mathcal{O}(\varepsilon^{-2})$, see
Theorem~\ref{thm:chimera-state}~[\ref{enum:thm-eps2-time}].
\item[{\it (Asymptotically stable)}] Chimera states, which persist for
all time, exist if the intra-star coupling configuration is sufficiently
sparse. (Theorem~\ref{thm:asymptotically_stable_chimera-state})
\end{description}

The choice of our specific model was inspired by earlier observations of Ko and Ermentrout \cite{Ko_2008_1,Ko_2008_2}, Vlasov \textit{et al.}~\cite{Vlasov2015} and Toenjes {\it et al.}~\cite{RalfPaper} that phase oscillators in star networks exhibit  coexisting coherent and incoherent states caused by shear-induced degree dependence of the oscillator frequency.

We obtain results on the persistence of mixed chimeras states obtained by weakly coupling two identical stars. In the uncoupled limit ($\varepsilon=0$), the system has an invariant manifold $M = M_+ \times M_-$, where $M_+$ denotes the invariant manifold associated with synchronous motion of the first star and $M_-$ the invariant manifold associated to the asynchronous attractor of the other star. We show that this invariant manifold is normally hyperbolic (Theorem~\ref{thm-NAIM}) which persists for sufficiently small intra-star coupling ($\varepsilon>0$), in the sense that the weakly coupled system has an invariant manifold  $M_\varepsilon$ that is close to $M$. The  chimera states persist for all time only if the corresponding trajectories of the system remain on $M_\varepsilon$ for all time, which may not be the case if $M_\varepsilon$ has a boundary $\partial M_{\varepsilon}$ through which trajectories can escape. We   may observe one of the following, see also Figure~\ref{fig:chimeras_birth} for a sketch:

\begin{figure}
\centering
\includegraphics[width=1.0\linewidth]{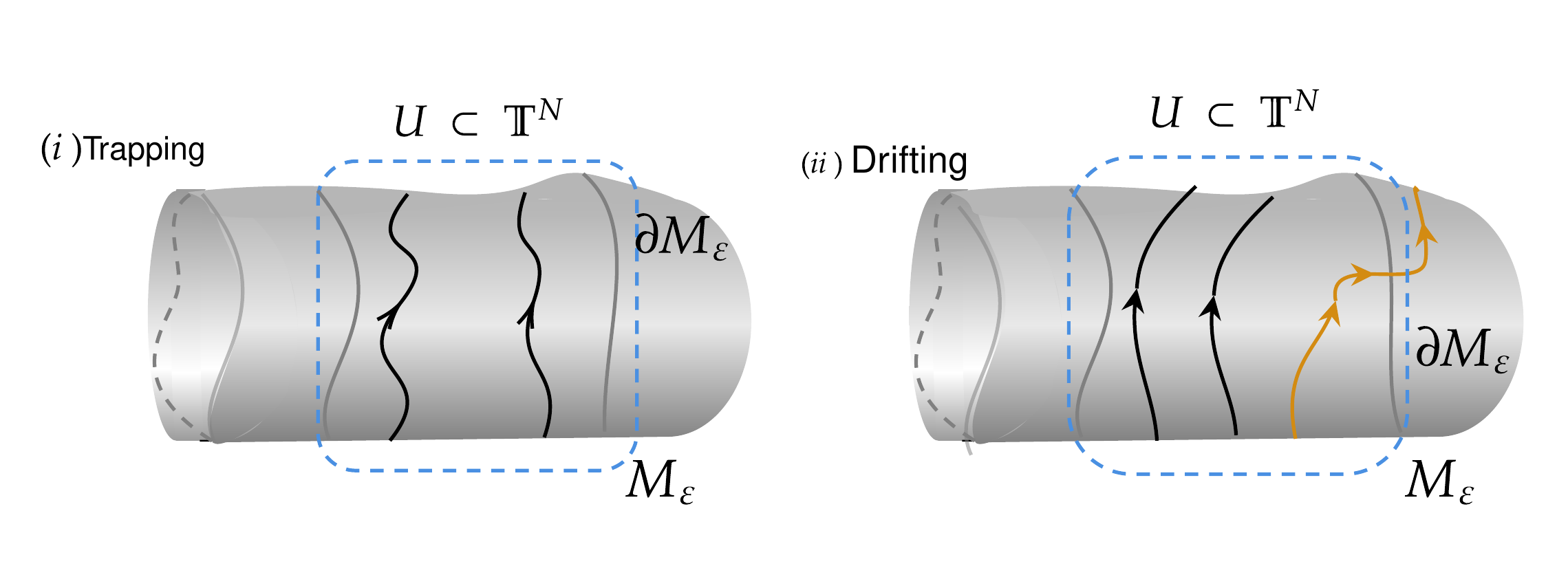}
\caption{The uncoupled stars have an invariant manifold product denoted by $M = M_s \times M_I$ with an invariant boundary $\partial M$ determined by $U \subset \mathbb{T}^{N}$. Once both stars are coupled under small coupling strength $\varepsilon$, the invariant manifold $M$ is perturbed to a nearby invariant manifold $M_{\varepsilon}$ with its boundary $\partial M_{\varepsilon}$. Simulations suggest two long-time behaviours: (i) trapping: the existence of a chimera state holds for large finite time; (ii) drifting: chimera  state ends up over finite time due to the dynamics leaves out the invariant manifold boundary $\partial M_{\varepsilon}$.}
\label{fig:chimeras_birth}
\end{figure}
\begin{itemize}
\item Stable chimera:  the trajectory cannot escape through the boundary $\partial M_{\varepsilon}$.
\item Metastable chimera:  the trajectory leaves the $M_{\varepsilon}$
in finite time by traversing $\partial M_{\varepsilon}$.
\end{itemize}
Our main results rely on estimates from normally hyperbolic invariant manifold and averaging theory to obtain lower bounds on the drift on $M_\varepsilon$. Before discussing the proofs of the main results in more detail, in the next section we first present the model in detail and present numerical observations that have motivated the theory developed in this paper.

\section{The coupled stars network}
\label{sec:coupled_stars_network}
In this section we describe the model and present numerical observations that motivated the main results obtained in this paper.
\subsection{Model}
\label{subsec:model}
Consider two coupled star networks labelled ``$+$" and ``$-$",
with hub states $\varphi_0^\pm\in S^1$  and leaf states
$\varphi_i^\pm \in S^1$, $i=1,\ldots, N$, where
$S^1\simeq \mathbb{R}/2\pi\mathbb{Z}$, with equations of motion
\begin{equation}\label{eq:coupled_stars_initial}
\begin{split}
\dot{\varphi}_0^+ &= \omega + \lambda \sum_{j=1}^{N}
H(\varphi_j^+, \varphi_0^+),\\
\dot{\varphi}_i^+ &= \omega + \lambda H(\varphi_0^+,
\varphi_i^+) + \epsilon  h(\varphi^+_i - \varphi^-_i),
\end{split}
\end{equation}

\noindent and likewise for the $``-"$ star. Here, $\omega$ is the natural frequency of the oscillators, $\lambda$ represents the corresponding coupling strength and $h$ is a diffuse pairwise coupling function between the two stars. We note that considering a general smooth $h$ for pairwise coupling will not affect the results. The function $H$, representing the coupling between the hub and leaves within each star, is the Kuramoto--Sakaguchi coupling  with shear
\begin{align}\label{eq:coupling_function_shear}
H(\varphi_j, \varphi_i) = c + \sin(\varphi_j - \varphi_i + \delta),
\end{align}
where $c \neq 0$ is the shear and $\delta \in (0, \pi/4)$ is the phase frustration. This coupling is related to the phase reduction description for a Hopf-Andronov bifurcation~\cite{Kuramoto1984}. Shear causes the frequency of the oscillators to depend on node degree.

As a measure of synchronization, we consider the order parameters
    \begin{equation}\label{eq:order_parameter}
        z^{+} = \frac{1}{N} \sum_{j=1}^N e^{i (\varphi_j^{+} - \varphi_0^{+})} \qquad z^{-} = \frac{1}{N} \sum_{j=1}^N e^{i (\varphi_j^{-} - \varphi_0^{-})},
    \end{equation}
 where $r^+ := |z^+|$ and $r^{-} := |z^{-}|$ are the absolute value.
\subsection{Change of coordinate and parametrization}
\label{sec:change_of_coordinate_and_parametrization}

We perform the variable substitution $\varphi_j \mapsto \tilde{\varphi}_j \coloneqq \varphi_j + \omega t$ for all $j = 0, \dots, N$
and a time rescaling $t \mapsto \tau \coloneqq \lambda c t$.
This rescales Eq.~\eqref{eq:coupled_stars_initial} and turns the
degree into the hub node's natural frequency, and gives
\begin{equation}\label{eq:explicit_relations}
\begin{split}
\dot{\varphi}_0^+ &= N + \frac{1}{c} \sum_{j=1}^{N}
\sin(\varphi_j^+ - \varphi_0^+ + \delta ),\\
\dot{\varphi}_i^+ &= 1 + \frac{1}{c} \sin(\varphi_0^+ -
\varphi_i^+ + \delta) + \frac{\epsilon}{\lambda c}  h(\varphi^+_i - \varphi^-_i),
\end{split}
\end{equation}
where we abuse notation by re-using variables $\varphi$ and denoting differentiation with respect to $\tau$ by $\cdot$ again. Let us denote $\sigma = 1/c$, $\varepsilon = \epsilon \sigma/\lambda$ and allow for a slight generalization, where $\beta>1$ is the frequency of the hub and $k_i = 1$ for $i = 1, \dots, N$. When $\beta = N$ we recover the previous Eq.~\eqref{eq:explicit_relations}.
Here $\beta$ works as parameter for the frequency
gap. So, the network  dynamics we analyze is given as
\begin{equation}\label{eq:coupled_stars}
\begin{split}
\dot{\varphi}_0^+ &= \beta + \beta \frac{\sigma}{N} \sum_{j=1}^{N}
\sin(\varphi_j^+ - \varphi_0^+ + \delta) \\
\dot{\varphi}_i^+ &= 1 + \sigma \sin(\varphi_0^+ - \varphi_i^+ + \delta)
+ \varepsilon  h(\varphi^+_i - \varphi^-_i), \quad i = 1,\dots, N,\\
\end{split}
\end{equation}
and likewise for the  ``$-$" star.

\subsection{Numerical observations}
\label{numerics}

We begin by reviewing the synchronization diagram for a single star and
then we proceed to numerically explore chimera states for two coupled stars.

\subsubsection{Isolated star network}

\begin{figure}
\centering
\includegraphics[width=0.7\linewidth]{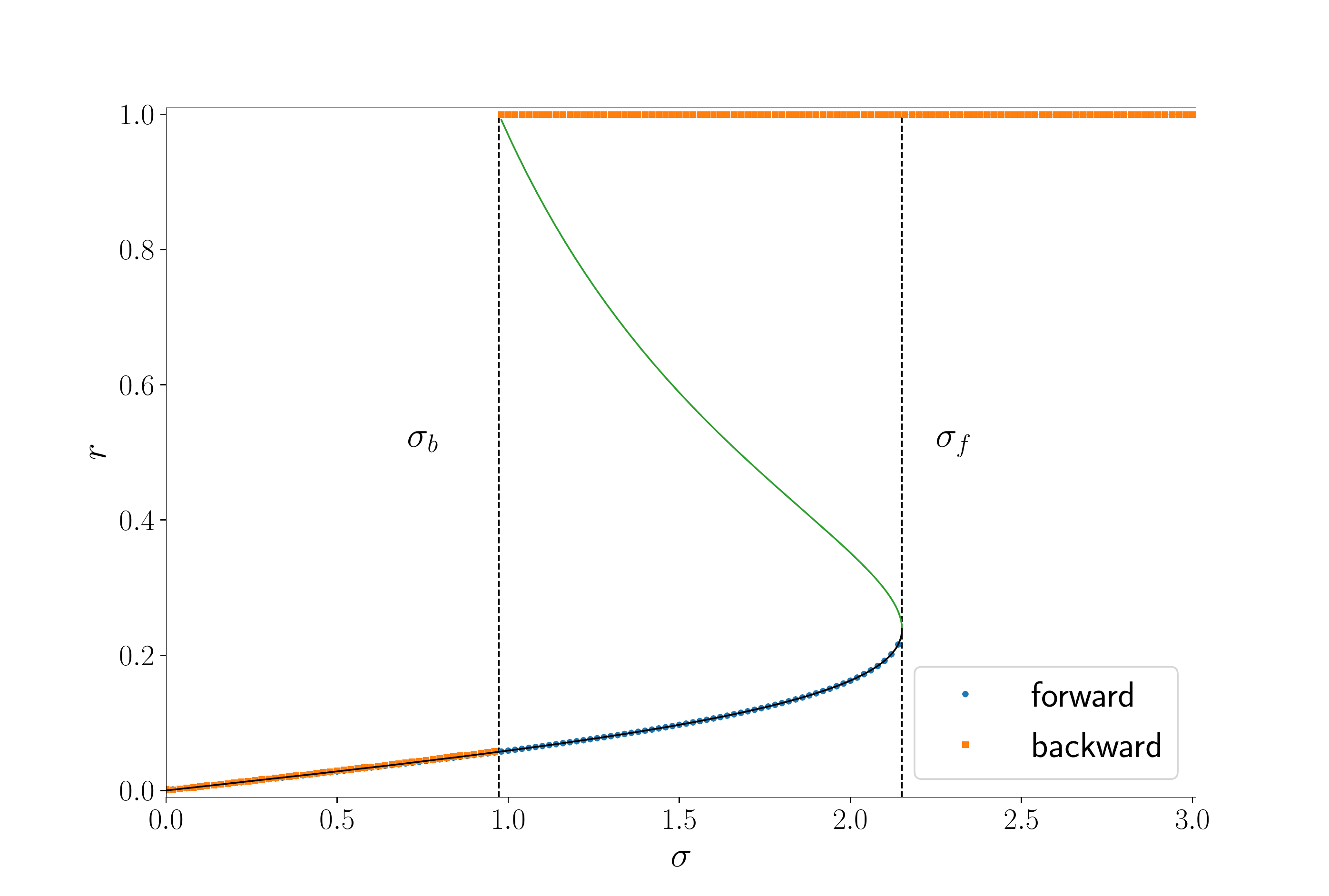}
\caption{Synchronization diagram for a single star with $N = 200$
leaves. The parameters are $\beta = 10$, $\delta = 0.3$ and at each step
the coupling strength was increased (decreased) by
$\Delta\sigma = 0.02$.  The region between the two vertical dashed
lines is where both coherent and incoherent states coexist. The green
solid line corresponds to the separatrix.}
\label{fig:single_star_diagram}
\end{figure}

We first obtain the bifurcation diagram of the order parameter as a
function of $\sigma$. Figure~\ref{fig:single_star_diagram} shows the
synchronization diagram for a star.

\textbf{Numerical procedure.} We integrate
Equations~\eqref{eq:coupled_stars} of one
star using an implicit Runge--Kutta method of order five (we use a solver described in~\cite{jitcxde}). Starting at $\sigma = 0$
with uniformly distributed initial phases in $(0, 2\pi)$ and evaluate the
order parameter $r$ in the stationary regime. Then we increase
adiabatically the coupling by $\Delta \sigma = 0.02$ and, using the
outcome of the last run as the initial condition, calculate the new value of
the stationary order parameter $r$ at $\sigma + \Delta \sigma$,
repeating these steps until a maximal value $\sigma$ is reached. This
curve is called the forward continuation. We can see that the order
parameters increases slowly and smoothly until
$\sigma = \sigma_{f}$, where it jumps discontinuously to $r = 1$.

The backward continuation is obtained by decreasing adiabatically the
coupling strength by steps of size $\Delta \sigma$ from the synchronous
state. We use the outcome of the last run as the initial condition for the
calculation of $r$ we add to the initial conditions a small random
uniform number drawn from the interval $(0,0.01)$. In the backward
continuation, the order parameter jumps back to the asynchronous branch.

\textbf{Main findings.} In Section~\ref{sec:sync_diagram}, we prove the
existence of such synchronization diagram via the theory of invariant
manifolds. Here we summarize the main points. In
Figure~\ref{fig:single_star_diagram}, the region between the dashed
lines corresponds to the values of the coupling critical couplings
\begin{align*}
\sigma_b = \frac{\beta-1}{1 + \beta \cos(2 \delta)} \quad
\mathrm{and} \quad \sigma_f =
\frac{\beta-1 }{\sqrt{1+2\beta\cos(2\delta)}}.
\end{align*}

The backward critical point $\sigma_b$ is calculated in
Section~\ref{sec:naim_sync}, where
we prove that it is a consequence of requiring the existence and
stability conditions for the synchronous state. The forward critical value
$\sigma_f$ is calculated in Section~\ref{sec:naim_async} analyzing
the differential equation of the M\"obius group $\alpha$ parameter.
This allows us to obtain the solid lines as well. The black solid line
corresponds to the $-$ sign and the green solid line to the $+$ in the
following equation
\begin{equation*}
\frac{(\beta-1)  \pm \sqrt{(\beta-1)^2 - \sigma^2(1+2\beta \cos(2\delta))}}{\sigma(1+2\beta \cos(2\delta))}.
\end{equation*}
Moreover, in Section~\ref{sec:sync_diagram_below_backward_critical_coupling}
we prove through the averaging principle that for small coupling strength
$\sigma < 1$ the incoherent state is a global attractor on the unit disk.

\subsubsection{Network of two symmetrically coupled stars}

We consider two stars coupled as in Figure~\ref{fig:coupled_stars}. Each
star has $N$ leaves, with the coupling strength between its leaves and its
hub given by $\sigma$. The coupling between them is mediated by
connecting, with strength $\varepsilon$, the $k$th leave of star $+$
to the $k$th leaf of star $-$.

The mechanism for generating a chimera state is a consequence of the
coexistence between synchronous and incoherent states. From identical
stars, we select a star to be at the synchronous state and another to be at
the incoherent. Both stars are coupled under an inter-coupling strength
$\varepsilon$ and function
\begin{equation*}
h(\varphi_j^{\pm}, \varphi_i^{\mp}) = \sin(\varphi_j^{\pm} - \varphi_i^{\mp}
+ \delta),~\quad ~i, j \in \{0, \dots N\}.
\end{equation*}
Starting from this prepared initial conditions, and from the arguments in
Section~\ref{sec:main-proof}, we show that the coupled system
presents similar macroscopic behavior as long as the \emph{perturbation}
caused by the coupling is small enough.

Examples of such chimera-like state are depicted in
Figure~\ref{fig:coupled_stars}. In Figure~\ref{fig:coupled_stars} where every leaf of the top star is connected to one of the bottom star,
with strength $\varepsilon = 0.007$. When the coupling between the stars is turned
on, the macroscopic behaviour persists where the order
parameter exhibits an oscillatory behavior with varying amplitude and
frequency.  We note that the oscillatory
behavior akin to the \emph{breathing chimera}~\cite{abrams2008}. \par 

%

\subsubsection{Scale free networks of coupled stars}

\begin{figure}
\centering
\includegraphics[width=0.6\linewidth]{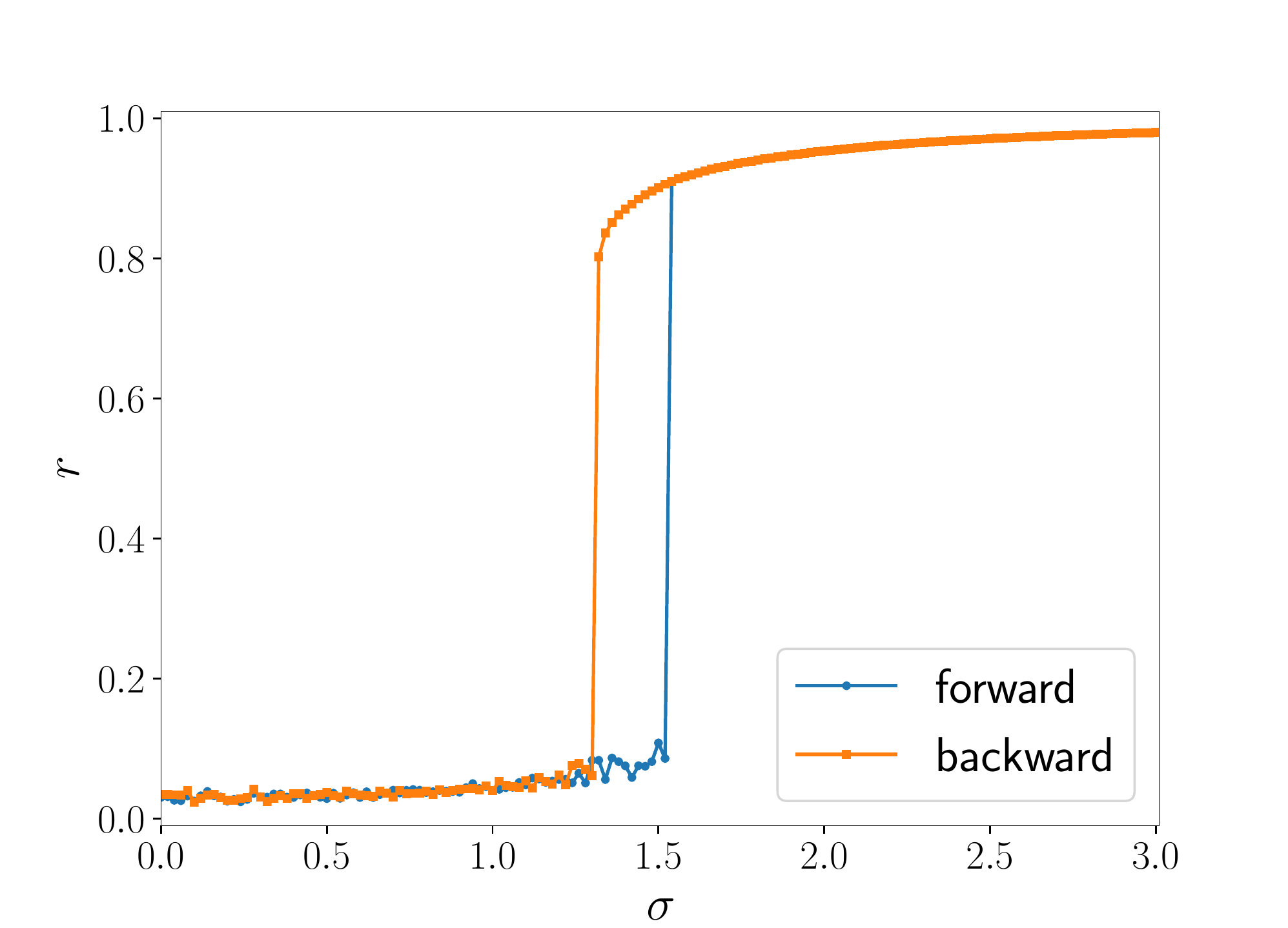}
\caption{Synchronization diagram for a single Barab\'asi--Albert
network with $N = 1000$ nodes, $\gamma \approx 2.5$,
$\langle k \rangle = 6$ and $\delta = 0.03$,
and at each step the coupling strength was increased (decreased) by
$\Delta\sigma = 0.02$. There is a coherent and incoherent states
coexistence region.}
\label{BA_sync_diagram}
\end{figure}

Star motifs are building blocks for complex networks with scale-free
networks \cite{barabasi1999}. Although our results for star graphs do
not apply to the scale-free network, we can use them as heuristics to
build chimera-like states in scale-free networks.  Consider the following
equations of motion on complex networks,

\begin{equation}\label{eq:kuramoto_complex_networks}
\dot{\varphi}_{i} = k_{i} + \sigma \sum_{j=1}^{N} A_{ij}
\sin \left(\varphi_{j} - \varphi_{i} + \delta \right), \quad i = 1, \dots, N,
\end{equation}

\noindent
where $\sigma$ is the coupling strength, $k_{i}$ is the degree of the
vertex $i$ and $A_{ij}$ is the adjacency matrix. The synchronization
diagram for a Barab\'asi--Albert~\cite{barabasi1999} network with
$1000$ vertices, degree distribution exponent $\gamma  \approx 2.5$
and mean degree $\langle k \rangle = 6$, performed with the same
method as was used in Figure~\ref{fig:single_star_diagram}
is shown in Figure~\ref{BA_sync_diagram}.

\begin{figure}
\centering
\includegraphics[width=1.0\linewidth]{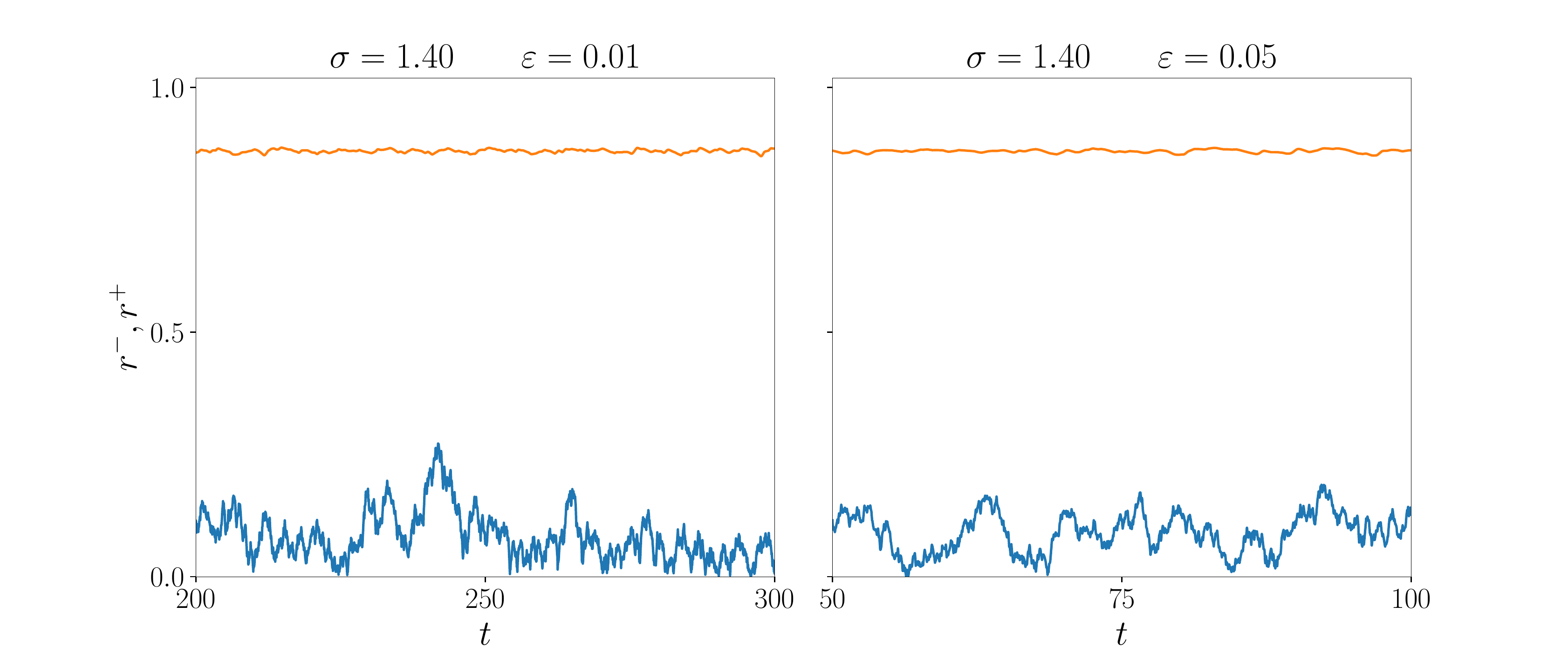}
\caption{ Time series for the order parameters $r^{+}$ and $r^{-}$
for two coupled Barab\'asi--Albert networks with $N = 1000$, degree
distribution exponent $\gamma \approx 2.5$ and mean degree
$\langle k \rangle = 6$. We indicate the intra-coupling strength $\sigma$
and inter-coupling strength $\varepsilon$. }
\label{coupled_BA_time_series}
\end{figure}

We coupled two Barab\'asi--Albert networks with the same degree
distribution exponent $\gamma$ and mean degree $\langle k \rangle$.
The $k$th vertex of one network is coupled to (and only to) to the $k$th
vertex of the other network. The coupling function used is the
Kuramoto--Sakaguchi coupling. Again $\sigma$ denotes the intra-coupling strength
and  the $\varepsilon$ inter-coupling.  Figure~\ref{coupled_BA_time_series}
shows chimera-like states in the two coupled Barab\'asi--Albert networks.
We chose initial conditions inside the hysteresis loop in
Figure~\ref{BA_sync_diagram}. The chimera breaking is observed
in the right column.

\section{Mathematical analysis}
\label{sec:mathematical_analysis_prelimanaries_and_results}
In the remainder of this paper, the main results, as informally discussed above in Section~\ref{sec:motivation_and_main_results} and summarized in Theorems~\ref{thm-NAIM}-\ref{thm:asymptotically_stable_chimera-state}, are proved. To this end we sketch the strategy to prove the main theorem.

\subsection{Strategy}
\label{sec:main_ideas_proof}

\begin{description}
\item[First Step:] {\it Constructing the Synchronous Manifold of a single star}.
There exists a manifold where
all phases are locked. In Section~\ref{sec:naim_sync} we prove that for certain
intervals of the parameters the manifold is normally attracting, see Proposition~\ref{prop:sync_state_NAIM}.

\item[Second Step:] {\it  Asynchronous manifold for a single star.} We
construct it as follows:

\subitem {\it Reduction in terms of M\"obius actions.} We rewrite the model
in terms of relative phases. In these coordinates,  dynamics are described
in a low-dimensional manner by the action of a three-parameter M\"obius
group, see Section~\ref{sec:ws_theory}. The Watanabe--Strogatz approach
provides a differential equation for the parameters of the
M\"obius group $\alpha \in \mathbb{C}$ and $\psi \in S^1$ keeping invariant a set of constants of motion.

\subitem {\it Invariant Manifold for asynchronous dynamics.} To construct the
manifold for the asynchronous dynamics we perform the following steps:

\subsubitem
(i) Construct open sets of initial conditions such that the
order parameter $z$ is $C^k$ close to  $\alpha$, see Section
\ref{sec:order_parameter_is_close_to_alpha} in Lemma~\ref{lem:z-alpha-close}. The differential equation for $\alpha$ depends on $z$ and by a
small perturbation we obtain an equation solely in terms of $\alpha$,
see Section~\ref{sec:closed_equation_alpha} in Theorem~\ref{thm:hyperbolicity_incoherent_branch}.
 By averaging theory, in Section~\ref{sec:sync_diagram_below_backward_critical_coupling}, we
prove that the perturbed dynamics of $\alpha$ enters in finite time a neighborhood of zero, see Theorem~\ref{thm:Asymptrho}. This determines that the bifurcation diagram of a single star network is in fact discontinuous for sufficiently small coupling strength.

\subsubitem (ii) By persistence of hyperbolic manifolds we obtain
the invariant manifold corresponding to the asynchronous dynamics
in the full equations, see Proposition~\ref{prop:async_state_NAIM}.

\noindent
\item[Third Step:]{\it Constructing Chimeras in coupled stars}. For the coupled
stars, the constants of motion may drift to the boundary of the
asynchronous manifold. When this happens we lose control over the
dynamics because the order parameter may no longer be in a neighborhood
of $\alpha$. So first, we prove the existence of chimera states in the couple stars system, see Theorem~\ref{thm-NAIM}. And finally, we show three distinct results about the time scale these chimera states exists depending on how generic we couple both stars.
\subitem{\it Short metastability.} We estimate an upper bound for how long it takes the perturbed dynamics of the constants of motion drift out the invariant manifold.
\subitem{\it Long metastability.} We assume the coupling function between both stars has sinusoidal form as an extra assumption. This implies the coupling term can be split up into two terms and restricted to the invariant manifold. This restriction introduces an effective coupling term of order $\varepsilon^2$ to the constants of motion dynamics.
\subitem{\it Asymptotically stable metastability.}
 We show that coupling a single leaf may generate a drift of the
conserved quantities, but this will happen only for a single constant.
That is, the drift does not cascade to other leaves. Hence, if only a
small fraction of leaves are coupled this drift is immaterial and the
order parameter will stay in a $C^k$ neighborhood of $\alpha$ and the
 chimera state is stable.
\end{description}

\subsection{Preliminaries}
\label{sec:preliminary_results}
In preparation for the subsequent analysis, we briefly recall some
established results from the theory normally hyperbolic invariant
manifolds and averaging.

\subsubsection{Normally attracting invariant manifolds (NAIM)}
\label{sec:NAIMS}

Our first ingredient  is the theory of normally
hyperbolic invariant manifolds~\cite{fenichel1971,HPS1977,eldering2013}.
Such manifolds can be viewed as generalizations of hyperbolic fixed
points and we use their property that they persist under small
perturbations of the dynamics. We shall only consider
normally \emph{attracting} invariant manifolds (NAIMs), i.e.\ those
which only have stable normal directions, not unstable ones.

Let $\dot{x} = f(x)$ with $x \in \R^n$ define a dynamical system; we
use $\R^n$ for simplicity of presentation, but it can be replaced by
a smooth manifold. A submanifold $M \subset \R^n$ is said to be a
normally attracting invariant manifold for $f$ if it is invariant
under the flow $\Phi^t$ of $f$ and furthermore the linearized flow
along $M$ contracts normal directions, and does so more strongly than
it contracts directions tangent to $M$. Our precise definition below
corresponds to that of ``eventual absolute normal hyperbolicity''
according to~\cite{HPS1977}.
\begin{definition}[Normally attracting invariant manifold]\label{def:NAIM}
Let $\dot{x} = f(x)$ with $x \in \R^n$ and $f \in C^1$ describe a
dynamical system with flow $\Phi\colon \R \times \R^n \to \R^n$.
Then a submanifold $M \subset \R^n$ is called a normally attracting
invariant manifold of $f$ if all of the following
conditions hold true:
\begin{enumerate}
\item $M$ is invariant, i.e.\ $\forall\; t \in \R\colon \Phi^t(M) = M$;
\item there exists a continuous splitting
\begin{equation}\label{eq:NAIM-split}
T_M \mathbb{R}^n = T M \oplus N
\end{equation}
of the tangent bundle $T \R^n$ over $M$ into the tangent
and a (stable) normal bundle, with continuous projections
$\pi_M,\,\pi_N$ and this splitting is invariant under the
tangent flow $T\Phi^t = T\Phi_M^t \oplus T\Phi_N^t$;
\item there exist real numbers $a < b \le 0$ and $C > 0$
such that the following exponential growth conditions hold
on the subbundles:
    \begin{eqnarray*}
      &\forall\; t\le 0,\,(m,v) \in T M\colon\quad
      \norm{T\Phi_M^t(m)\,v} \le C\,e^{b\,t}\,\norm{v},\\
      &\forall\; t\ge 0 ,\,(m,v) \in N \colon\qquad\!
      \norm{T\Phi_N^t(m)\,v} \le C\,e^{a\,t}\,\norm{v}.
    \end{eqnarray*}
  \end{enumerate}
\end{definition}
Note that if we have two dynamical systems $f_1,f_2$ which each have a
NAIM $M_1,M_2$ respectively and $\max(a_1,a_2) < \min(b_1,b_2)$, then
$M_1 \times M_2$ is a NAIM for the product system. In our case, the
dynamics on the invariant manifolds will be (close to) neutral, so
$b_{1,2} \approx 0$ and hence the product system has a NAIM too.

A main result on normal hyperbolicity is that NAIMs persist under
small perturbations of the vector field $f$. That is, we have the
following theorem, see~\cite[Thm.~1]{fenichel1971}
or~\cite[Thm.~4.1]{HPS1977}.
\begin{theorem}[persistence of NAIMs]\label{thm:persist-NAIM}
  Let $M \subset \R^n$ be a compact normally attracting invariant
  manifold for $f$. Then there exists an $\varepsilon > 0$ such that
  for any vector field $\tilde{f}$ with
  $\norm{\tilde{f} - f}_{C^1} \le \varepsilon$, there exists a unique
  manifold $\tilde{M}$ that is diffeomorphic and
  $\mathcal{O}\bigl(\norm{\tilde{f} - f}_{C^1}\bigr)$-close to $M$ and
  invariant under $\tilde{f}$. Furthermore, $\tilde{M}$ is a NAIM for
  $\tilde{f}$.
\end{theorem}

The closeness of $\tilde{M}$ to $M$ can in general be expressed
as follows: $\tilde{M}$ can be described by the graph of a section
of the normal bundle of $M$, and this section is $C^1$-small. In
our case $M$ can explicitly be given by the graph of a
function, hence $\tilde{M}$ will be given by the graph of a
$C^1$-small perturbation of this function.

In our case the NAIM $M$ will have a boundary.
Theorem~\ref{thm:persist-NAIM} still holds true as long as the
vector field $f$ is pointing strictly outward at the boundary $\partial M$
(also called `overflowing invariant'), see~\cite{fenichel1971}. In our
case $M$ will not be overflowing invariant. This can be overcome
by artificially modifying $f$ near the boundary $\partial M$. A slightly
simpler approach uses~\cite[Thm.~4.8]{eldering2013}: we modify $f$
such that it preserves a manifold $S$ that transversely intersects
$M$ with $S \cap M = \partial M$. For example, $S$ could be the
local normal bundle of $M$ restricted to $\partial M$.

\subsubsection{Averaging principle}
\label{sec:averaging_principle}

The averaging principle will play role in our analysis. Consider
the \textit{original} system given by
\begin{align}\label{eq:original_system_thm}
\begin{array}{lcl}
\dot{x} & = & \varepsilon F(x, \eta), \\
\dot{\eta} & = & \Omega(x, \eta) + \varepsilon~G(x, \eta), \quad
(x(0), \eta(0)) = (x_0, \eta_0), x \in D,  \eta \in S^1.
\end{array}
\end{align}
and the \textit{averaged} system given by
\begin{align}\label{eq:averaged_system_ode}
\begin{array}{lcl}
\dot{y} & = & \varepsilon F^1(y), \\
\dot{\zeta} & = & \Omega(y, \zeta), \quad (y(0), \zeta(0))
= (x_0, \zeta_0) \in U \times S^1,
\end{array}
\end{align}
where
\begin{equation*}
F^1 (y) = \frac{\frac{1}{2 \pi} \int_0^{2\pi}
\frac{F(y, \vartheta)}{\Omega(y, \vartheta)} d\vartheta }
{\frac{1}{2 \pi} \int_0^{2 \pi} \frac{d \vartheta}{\Omega(y,\vartheta)}}
\end{equation*}
is the first order averaged vector field. If the averaged system has
an asymptotically stable fixed point, solutions of the original system
are approximated by averaged ones. That is, we have the following
adapted version of the theorem, see~\cite[5, Thm.~7.2]{Chicone2006}
or~\cite[7, Lemma.~7.10.1-2, Thm.7.10.4]{Sanders2007Averaging}.

\begin{theorem}\label{thm:averaging_attractor}\emph{[Averaging
theorem with attraction ]}
Let $D\subset \mathbb{R}^n$ be an open set and consider the
equation~\eqref{eq:original_system_thm}. Suppose that the averaged
system Equation~\eqref{eq:averaged_system_ode} has an asymptotically
stable fixed point $y = 0$, $F^1 \in C^{1}(D)$ and has a domain of
attraction $D_0 \subset D$. For any compact $K\subset D_0$ there
exists a $\varepsilon_0$ and $c > 0$ such that for all $x_0 \in K$ and
for each $\varepsilon < \varepsilon_0$
\begin{equation*}
  \|x(t) - y(t)\| \leq c\varepsilon, \quad 0 \leq t < \infty.
\end{equation*}
\end{theorem}

{Bifurcation diagram of an isolated star network}
\label{sec:sync_diagram}

\subsection{A NAIM for the synchronous state}
\label{sec:naim_sync}
The synchronous state is characterized by $r = 1$. The phase locking
manifold
\begin{equation*}
M_{\phi} = \{\phi_1 = \ldots = \phi_N = \phi \in S^1, \phi - \varphi_0
= \Delta \phi \in S^1\} \subset \mathbb{T}^{N + 1}
\end{equation*}
satisfies this condition on the order parameter. The stability of $M_{\phi}$ depends on the frequency mismatch, the phase frustration $\delta$ and
intra-coupling strength $\sigma$~\cite{Vlasov_kuramoto_japanese_drums_2015}.
The following hypothesis is necessary to prove existence of the
coherent manifold in
Proposition~\ref{prop:sync_state_NAIM}.
\begin{hypothesis}\label{hyp:stable-state-cond}
Assume that
\begin{equation}\label{eq:stable-sync-cond}
\sigma > \sigma_b := \frac{\beta-1}{1 + \beta \cos(2 \delta)}
\qquad\textrm{and}\qquad
\textrm{arg}(v) - \sin^{-1}\Bigl(\frac{1-\beta}{\sigma\norm{v}}\Bigr)
\in \Bigl( -\frac{\pi}{2} , \frac{\pi}{2} \Bigr),
\end{equation}
where $v = (\beta \sin (2 \delta), 1 + \beta \cos (2 \delta)) \in \R^2$,
and hence $\norm{v} = \sqrt{1 + \beta^2 + 2\beta\cos(2\delta)}$.
\end{hypothesis}

\begin{proposition}\label{prop:sync_state_NAIM}
  Let $\sigma > 0$, $\beta > 1$ and $\delta \in (0,\frac{\pi}{4})$. The synchronized state
  \begin{equation}\label{eq:invar-coherent}
    M_C = \{\, \phi_1 = \ldots = \phi_N = \phi^C \in S^1 \,\}
  \end{equation}
  is a normally attracting invariant manifold for the
  dynamics~\eqref{eq:single_star_rel_phase} if and only if
  Hypothesis~\ref{hyp:stable-state-cond} holds and when $\phi^C$ is
  given by
  \begin{equation}\label{eq:stable-sync-phi}
    \phi^C = \delta - \pi + \arctan\left(\frac{1 + \beta \cos(2 \delta)}{\beta
    \sin(2 \delta)} \right) + \arccos \left(\frac{\beta - 1}{\sigma \norm{v}}\right).
  \end{equation}
\end{proposition}

\begin{figure}[htb]
  \centering
  \input{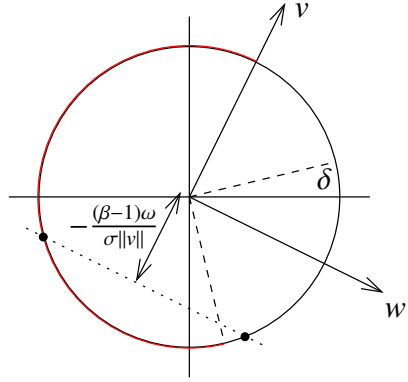}
  \caption{A geometric representation of the existence and stability
    of the synchronized state. The dots indicate the fixed points for
    $\phi^C \in S^1$ and the red arc is the unstable region. The
    vectors $v$ and $w$ and the two dashed lines are perpendicular.}
  \label{fig:geom-stable-sync}
\end{figure}

\begin{proof}
  First we search for the value of $\phi = \phi^C$ such that $M_C$ is
  an invariant manifold. Evaluating~\eqref{eq:single_star_rel_phase}
  on $M_C$ we see that this amounts to
  \begin{equation*}
    \dot{\phi} = 1 - \beta - \sigma \sin(\phi + \delta)
                             - \beta \sigma \sin(\phi - \delta) = 0.
  \end{equation*}
  This can be rewritten as
  \begin{equation}\label{eq:fixed-sync-state}
    \inner{v}{\bigl(\cos(\phi^{C} - \delta),\sin(\phi^{C} - \delta)\bigr)}
    = -\frac{\beta-1 }{\sigma}
  \end{equation}
  with $v = (\beta \sin (2 \delta), 1 + \beta \cos (2 \delta))$.
  Hence we have solutions if
  \begin{equation}\label{eq:sync-state-ineq}
    \norm{v} = \sqrt{1 + \beta^2 + 2\beta\cos(2\delta)}
    \ge \frac{\beta-1 }{\sigma}.
  \end{equation}

  To evaluate stability of $M_C$, let
  $(\dot{\phi}_1,\ldots,\dot{\phi}_N,\dot{\varphi}_H) =
  F(\phi_1,\ldots,\phi_N,\varphi_0)$
  denote the system~\eqref{eq:single_star_rel_phase}. Then we have
  the total derivative
  \begin{equation*}
    DF|_{M_C} = -\sigma \cos(\phi^C+\delta)
    \left(\begin{array}{cccc}
      1 & & & \emptyset \\
        & \ddots &   & \\
        &        & 1 & \\
     \emptyset & &   & 0
   \end{array}\right)
                -\sigma \frac{\beta}{N} \cos(\phi^C-\delta)
    \left(\begin{array}{cccc}
      1 & \cdots &  1 & 0 \\
      \vdots & \ddots & \vdots & \vdots \\
      1 & \cdots &  1 & 0 \\
     -1 & \cdots & -1 & 0
   \end{array}\right).
  \end{equation*}
  Thus, we find that $(0, \ldots, 0, 1)$ is an eigenvector with
  eigenvalue $0$, that is, the direction along $M_C$ is neutral. In
  the directions transversal to $M_C$ we have the eigenvector
  \begin{equation*}
    v_1 = \Bigl(1,\ldots,1,
     -\frac{\beta\cos(\phi^C+\delta)}{\cos(\phi^C-\delta)+
     \beta\cos(\phi^C+\delta)}\Bigr)
   \end{equation*}
  with eigenvalue $\lambda_1 = -\sigma\bigl[\cos(\phi^C-\delta)
  + \beta\cos(\phi^C+\delta)\bigr]$ and $N-1$ independent eigenvectors
  $v_i = (0\ldots, 1, -1, 0\ldots, 0)$ with eigenvalue
  $\lambda_2 = -\sigma\cos(\phi^C-\delta)$. So $M_C$ is normally
  attracting if $\sigma > 0$ and
  \begin{equation*}
    \cos(\phi^C-\delta) > 0
    \qquad\textrm{and}\qquad
    \cos(\phi^C-\delta) + \beta\cos(\phi^C+\delta) > 0.
  \end{equation*}
  These conditions can be rewritten as
  \begin{equation}\label{eq:stable-sync-state}
    \begin{split}
          \inner{w}{\bigl(\cos(\phi^C - \delta),\sin(\phi^C - \delta )\bigr)} &> 0,
    \end{split}
  \end{equation}
  with $w = (1 + \beta \cos (2 \delta), -\beta \sin (2 \delta))$
  equal to $v$ rotated clockwise over 90 degrees.
  Interpreting the conditions~\eqref{eq:fixed-sync-state}
  and~\eqref{eq:stable-sync-state} geometrically with $\phi^C \in S^1$,
  see Figure~\ref{fig:geom-stable-sync}, we see that we have at most
  two synchronized states. Since $v \perp w$, precisely one satisfies
  the second stability condition of~\eqref{eq:stable-sync-state} if
  and only if~\eqref{eq:sync-state-ineq} holds as a strict
  inequality. Then that state $M_C$ is actually stable if $\cos(\phi^C-\delta)>0$
  holds, which can be expressed as the second condition
  in~\eqref{eq:stable-sync-cond} and solving for $\phi^C$
  yields~\eqref{eq:stable-sync-phi}.

  The backward critical coupling $\sigma_b$ is obtained solving
  equation $\cos(\phi^C - \delta)$ = 0 with $\phi^C$ expression
  obtained. This yields the first condition at~\eqref{eq:stable-sync-cond}.
\end{proof}

\subsection{A NAIM for the asynchronous state}\label{sec:naim_async}

Again let us consider a single star network ($\varepsilon = 0$) in
system~\eqref{eq:coupled_stars}. We change to a new set of
coordinates defined by $\phi_k = \varphi_k - \varphi_0$. The
$\phi_k$ measure the phase difference between the leaves and
the hub, and the resulting equations are
\begin{equation}\label{eq:single_star_rel_phase}
  \begin{split}
    \dot{\phi}_i    &= 1 - \beta - \sigma \sin(\phi_i - \delta)
-\beta \frac{\sigma}{N} \sum_{j=1}^{N} \sin(\phi_j + \delta), \\
    \dot{\varphi}_0 &= \beta + \beta \frac{\sigma}{N}
    \sum_{j=1}^{N} \sin(\phi_j + \delta).
  \end{split}
\end{equation}
In the original coordinates the system had a global $S^1$-symmetry
by which we have effectively reduced the system, since the equations
for the relative angles $\phi_i$ are now decoupled from $\varphi_0$.

To find a stable asynchronous state, we shall change the equations
for a single star from relative phase form~\eqref{eq:single_star_rel_phase}
to new variables introduced by
 Watanabe--Strogatz~\cite{watanabe1994}.

\subsubsection{Low-dimension description}
\label{sec:ws_theory}

Watanabe and Strogatz~\cite{watanabe1994}, showed
that a large class of systems of $N$ identical coupled oscillators can
be described by a coupled set of only three ODEs. They achieved this
through a time dependent change of variables, that explicitly shows
that the system possesses $N - 3$ constants of motion.
In~\cite{marvel2009} Marvel, Strogatz and Mirollo provided an
explanation in terms of M\"obius actions. We follow~\cite{marvel2009}
and consider systems of $N$ identical coupled phase oscillators
with equations of motion that can be put in the form
\begin{equation}\label{eq:phase_ws_form}
  \dot{\phi}_j = f e^{i\phi_j} + g + \bar{f} e^{-i\phi_j}.
\end{equation}
Here $f$ is a smooth complex-valued function of the phases
$\phi_{1}, \phi_{2}, \ldots, \phi_{N} \in S^1$, $\bar{f}$ is its complex
conjugate, and $g$ is in its turn a real-valued function of the phases
$\phi_{1}, \phi_{2}, \ldots, \phi_{N} \in S^1$. The phases $\phi_{j}$
evolve in time according to the action of the M\"obius\footnote{The
complete M\"obius group consists of all fractional linear transformations
$M(w) = \frac{aw+b}{cw+d}$ of the complex plane such that
$ad-bc \neq 0$ but we only consider the subgroup that preserves
the unit circle. As in related literature we shall refer to this subgroup
as the M\"obius group.} group,
\begin{equation}\label{eq:mobius_action}
  e^{i\phi_{j}(t)} = M_{\alpha(t),\psi(t)}(e^{i\theta_j}),
\end{equation}
where $\theta_j$ are a set of constant (time independent) angles,
and we use the group parametrization $M_{\alpha,\psi}$ and its
action on $w \in \C$ given by
\begin{equation}\label{eq:mobius_group}
  M_{\alpha,\psi}(w) = \frac{\alpha + e^{i\psi} w}{1 + \bar{\alpha} e^{i\psi} w},
\end{equation}
with $\psi \in S^1$ and $\alpha \in \D = \{z \in \C \mid \abs{z} < 1 \}$.

The time evolution of the phases in~\eqref{eq:mobius_action} satisfies
the equations~\eqref{eq:phase_ws_form} as long as $\alpha(t)$ and
$\psi(t)$ are solutions of the differential equations
\begin{equation}\label{eq:mobius_evolution}
  \begin{split}
  \dot{\alpha} &= i \left( f\alpha^2 + g\alpha + \bar{f} \right),\\
  \dot{\psi}   &= f\alpha + g + \bar{f}\bar{\alpha}.
  \end{split}
\end{equation}
The functions $f$ and $g$ still implicitly depend on
$\phi_{1}, \phi_{2}, \ldots, \phi_{N} \in S^1$, so we have to
use~\eqref{eq:mobius_action} to express the $\phi_j$'s in terms of
$\alpha,\psi$, and the $\theta_j$'s and obtain a closed system of
equations for $\alpha,\psi$. For suitable models, we can express $f$
and $g$ in terms of the order parameter $z$, i.e.\ $f = f(z,\bar{z})$
and $g = g(z,\bar{z})$, where $z$ is given by
\begin{equation}\label{eq:order-param-func}
  z(\alpha,\psi,\theta)
  := \frac{1}{N} \sum_{j=1}^N e^{i \phi_j}
   = \frac{1}{N} \sum_{j=1}^N M_{\alpha,\psi}\bigl(e^{i \theta_j}\bigr).
\end{equation}

We now reformulate the Watanabe--Strogatz change of variables in a
more formal geometric setting. Let $\dot{\phi} = X(\phi)$ with
$X \in \mathfrak{X}(\T^N)$ denote the vector field
of~\eqref{eq:phase_ws_form} and let
$(\dot{\alpha},\dot{\psi},\dot{\theta}) = \hat{X}(\alpha,\psi,\theta)$
with $\hat{X} \in \mathfrak{X}(Q)$ denote the vector field associated
to the Watanabe--Strogatz coordinates
$(\alpha,\psi,\theta) \in Q := \D \times S^1 \times \T^N$, i.e.\ we
have $\dot{\theta} = 0$ and $\dot{\alpha},\dot{\psi}$ given
by~\eqref{eq:mobius_evolution}. Let $\Phi^t$ and $\hat{\Phi}^t$ denote
the associated flows. Then we have the commuting diagram
\begin{equation}\label{eq:diagram-WS-coords}
  \xymatrix{
    Q \ar[d]^\pi \ar[r]^{\hat{\Phi}^t} & Q \ar[d]^\pi \\
    \T^N \ar[r]^{\Phi^t} & \T^N
  }
\end{equation}
where $\pi\colon Q \to \T^N$ is the submersion
$\phi = \pi(\alpha,\psi,\theta)$ defined by~\eqref{eq:mobius_action}.
Since $\pi(\alpha,\psi,\,\cdot\,)\colon \T^N \to \T^N$ is given by the
diagonal action of $M_{\alpha,\psi}$, it is a diffeomorphism and hence
$\pi$ indeed a submersion. Note that although we loosely speak of
$\pi$ as a coordinate transformation, strictly speaking it is not,
since it is not injective. The diagram~\eqref{eq:diagram-WS-coords}
implies that $T\pi \circ \hat{X} = X \circ \pi$ and hence
$X = T\pi \circ \hat{X} \circ \pi^{-1}$ for any right-inverse
$\pi^{-1}$. On the other hand, a vector field $Y\in\mathfrak{X}(\T^N)$
can be lifted to a vector field
$\hat{Y} = R \circ X \circ \pi \in \mathfrak{X}(Q)$
satisfying~\eqref{eq:diagram-WS-coords}, given a choice of a
right-inverse $R$ of $T\pi$. We now choose
\begin{equation}\label{eq:right-inverse-WS}
  R(\alpha,\psi,\theta) = T_\theta \bigl[\pi(\alpha,\psi,\,\cdot\,)^{-1}\bigr]\colon
  T_\phi \T^N \to T_\theta \T^N \subset T_{(\alpha,\psi,\theta)} Q
\end{equation}
with $\phi = \pi(\alpha,\psi,\theta)$. Let $M_{\alpha,\psi}'$ denote
the derivative of the action on $S^1$, then we have
\begin{equation}\label{eq:right-inverse-explicit}
  R(\alpha,\psi,\theta) =
  \operatorname{diag}\bigl(M_{\alpha,\psi}'(\theta_1)^{-1},\ldots,M_{\alpha,\psi}'(\theta_N)^{-1}\bigr).
\end{equation}
The upshot of this is that if we perturb the original vector field $X$
to $X + \varepsilon Y$, for example when introducing a coupling
between the stars, then we can lift $Y$ to
$\hat{Y} \in \mathfrak{X}(Q)$ such that $\hat{Y}$ is given by
\begin{equation}\label{eq:lifted-vectorfield}
  \dot{\alpha} = 0, \qquad
  \dot{\psi}   = 0, \qquad
  \dot{\theta}_j = M_{\alpha,\psi}'(\theta_j)^{-1} Y_j(\pi(\alpha,\psi,\theta)).
\end{equation}
In the case that only some of the components $Y_j$ are nonzero, it
follows that only the associated $\theta_j$'s are not conserved
anymore by $\hat{Y}$, while the other $\theta_j$'s are still
conserved.

Note that formalizing the Watanabe--Strogatz method itself this way,
leads to a different lifted vector field; indeed, with our choice of
lift we would obtain $\dot{\alpha} = 0$ and $\dot{\psi} = 0$ and gain
nothing. Instead, if $X \in \mathfrak{X}(\T^N)$ is given by
\begin{equation}\label{eq:WS-lift}
  \dot{\phi}_j = X_j(\phi) = g(\phi) + f(\phi) e^{i\phi_j} + \overline{f(\phi)} e^{-i\phi_j},
\end{equation}
then the Watanabe--Strogatz approach tells us that
$\hat{X} \in \mathfrak{X}(Q)$ given by~\eqref{eq:mobius_evolution} and
$\dot{\theta} = 0$ is a lift of $X$. If $X$ depends on extra
parameters then this directly translates to the same parameter
dependence for $\hat{X}$. The upshot is that we can apply different
lifts to parts of a vector field $X$, thereby obtaining the most
convenient choice of lifted vector field $\hat{X}$.

\subsubsection{Complex order parameter is approximated by alpha}
\label{sec:order_parameter_is_close_to_alpha}
We write the sine functions in~\eqref{eq:single_star_rel_phase} in exponential form.
This shows that we can recast the equations in the Watanabe--Strogatz form~\eqref{eq:phase_ws_form} with $f$ and $g$ given by
\begin{equation}\label{eq:WS-functions-fg}
  \begin{split}
    f(z,\bar{z}) &= \frac{-\sigma e^{-i\delta}}{2i}, \\
    g(z,\bar{z}) &= 1 - \beta  - \sigma \beta
                    \left( \frac{z e^{i\delta} - \bar{z} e^{-i\delta}}{2i} \right).
  \end{split}
\end{equation}

For a more detailed analysis of a stable incoherent state in the system,
we shall use that $z \approx \alpha$ when the $\theta$'s are close
to being uniformly distributed on $S^1$ (called `splay states'), and
when we are close to the thermodynamic limit, i.e.\ when the number
of leaves $N$ is sufficiently large. Let us make this more precise.
\begin{definition}\label{def:unif-distrib}
  We say that $\theta \in \T^N$ is uniformly distributed if
  there exists a $\varphi \in S^1$ such that
  $\theta_j = \frac{2 \pi j}{N} + \varphi$ up to permutation. We
  denote the set of all such $\theta$ by $\Theta$.
\end{definition}

\begin{lemma}\label{lem:z-alpha-close}
  Fix $0 < r < 1$ and consider the disc
  $\D_r = \{ \alpha \in \C \mid \abs{\alpha} < r \}$. For each
  $k \ge 1$ and $\varepsilon > 0$ there exists for all $N$
  sufficiently large an open neighborhood $U$ of\/
  $\Theta \subset \T^N$ such that the map
  $\D_r \to \C,\; \alpha \mapsto z(\alpha,\psi,\theta) - \alpha$, with
  $z(\alpha,\psi,\theta)$ given by~\eqref{eq:order-param-func},
  is smaller than $\varepsilon$ in $C^k$-norm, uniformly for all
  $(\psi,\theta) \in S^1 \times U$.
\end{lemma}

\begin{proof}
  Consider an element $\theta \in \Theta$ and wlog.\ assume that
  it is of the form as in Definition~\ref{def:unif-distrib} without
  permutation of the indices $j$ and that these run from $0$ to $N-1$.

  Below, we use the geometric series expansion
  $(1 + x)^{-1} = \sum_{l \ge 0} (-x)^l$ with
  $x = \bar{\alpha} e^{i \psi} e^{i \theta_j}$, which converges for
  any $\abs{\alpha} < 1$. We find
  \begin{eqnarray*}
    z(\alpha,\psi,\theta)
    &= \frac{1}{N} \sum_{j=0}^{N-1} \frac{\alpha + e^{i \psi} e^{i \theta_j}}{1 +
    \bar{\alpha} e^{i \psi} e^{i \theta_j}} \\
    &= \frac{1}{N} \sum_{j=0}^{N-1} \alpha \Bigl(1 + \alpha^{-1}
    e^{i(\frac{2 \pi j}{N} + \varphi + \psi)}\Bigr)
                   \sum_{l \ge 0} \Bigl(-\bar{\alpha} e^{i\bigl(\frac{2 \pi j}{N} +
                   \varphi + \psi\bigr)}\Bigr)^l \\
    &= \alpha \frac{1}{N} \sum_{j=0}^{N-1} \Bigl(1 + \alpha^{-1}
    e^{i(\frac{2 \pi j}{N} + \varphi + \psi)}\Bigr)
      \Biggl(\sum_{l=0}^{N-2} (-\bar{\alpha})^l e^{i\bigl(\frac{2 \pi j l}{N} +
       l(\varphi + \psi)\bigr)} + R_{N-2,j}(\alpha) \Biggr) \\
    &= \alpha \frac{1}{N} \sum_{l=0}^{N-2} \sum_{j=0}^{N-1}
          \Bigl[(-\bar{\alpha})^l e^{i\bigl(\frac{2 \pi j l}{N} + l(\varphi + \psi)\bigr)}
              + \frac{(-\bar{\alpha})^l}{\alpha} e^{i\bigl(\frac{2 \pi j (l+1)}{N} +
              (l+1)(\varphi + \psi)\bigr)} \Bigr] \\
    & \hspace{1cm}
      + \frac{\alpha}{N} \sum_{j=0}^{N-1} \Bigl(1 + \alpha^{-1}
      e^{i(\frac{2 \pi j}{N} + \varphi + \psi)}\Bigr) R_{N-2,j}(\alpha) \\
    &= \alpha \frac{1}{N} \sum_{l=0}^{N-2} \delta_{l,0} \sum_{j=0}^{N-1}
    (-\bar{\alpha})^l e^{i\bigl(\frac{2 \pi j l}{N} + l(\varphi + \psi)\bigr)}
      + \frac{\alpha}{N} \sum_{j=0}^{N-1} \Bigl(1 + \alpha^{-1}
      e^{i(\frac{2 \pi j}{N} + \varphi + \psi)}\Bigr) R_{N-2,j}(\alpha) \\
    &= \alpha + \frac{1}{N} \sum_{j=0}^{N-1} \Bigl(\alpha
    + e^{i(\frac{2 \pi j}{N} + \varphi + \psi)}\Bigr) R_{N-2,j}(\alpha).
  \end{eqnarray*}
  To obtain the stated result we have to bound the second term in
  $C^k$-norm. We have
  \begin{equation*}
    R_{N-2,j}(\alpha)
    = \sum_{l \ge N-1} \Bigl(-\bar{\alpha} e^{i\bigl(\frac{2 \pi j}{N}
    + \varphi + \psi\bigr)}\Bigr)^l
    = (-\bar{\alpha})^{N-1} \sum_{l \ge 0} (-\bar{\alpha})^l
    \Bigl(e^{i\bigl(\frac{2 \pi j}{N} + \varphi + \psi\bigr)}\Bigr)^{l+N-1}
  \end{equation*}
  and note that the sum still defines a function with radius of
  convergence one, so this function and its derivatives up to order
  $k$ are uniformly bounded for $\alpha \in \D_r$ by some number
  $B_k > 0$, hence $\norm{R_{N-2,j}}_{C^k} \le k r^{N-1} B_k$. Since the
  function $F_j(\alpha) = \alpha + e^{i(\frac{2 \pi j}{N} + \varphi + \psi)}$
  multiplying $R_{N-2}(\alpha)$ has $C^1$-norm bounded by $2 + r$ and
  all higher derivatives zero, we find
  \begin{eqnarray*}
    \norm[\Big]{\frac{1}{N} \sum_{j=0}^{N-1} F_j \cdot R_{N-2,j}}_{C^k}
    &\le \max_{j \in [0,N-1]} \sum_{m=0}^k \Bigl(
            \norm{F_j}_{C^0} \norm{R_{N-2,j}}_{C^k}
        + m \norm{F_j}_{C^1} \norm{R_{N-2,j}}_{C^{k-1}} \Bigr) \\
    &\le \Bigl(k + \frac{k(k+1)}{2}\Bigr) \max_{j \in [0,N-1]}
    \norm{R_{N-2,j}}_{C^k} \\
    &\le \frac{k^2(k+3)}{2} (2 + r) r^{N-1} B_k.
  \end{eqnarray*}
  Thus, for fixed $r < 1$ and $k \ge 1$ and given $\varepsilon > 0$,
  we see that there exists an $N_0 > 0$ such that for all $N \ge N_0$
  \begin{equation*}
    \norm{\alpha \mapsto z(\alpha,\psi,\theta) - \alpha}_{C^k}
    < \frac{\varepsilon}{2}
    \qquad\textrm{on } \D_r \times S^1 \times \Theta.
  \end{equation*}
  Next, since this function is $C^\infty$ (and actually analytic
  w.r.t.~$\alpha,\bar{\alpha},\psi,\theta$) it follows that there
  exists an open neighborhood $U \supset \Theta$ such that
  \begin{equation*}
    \norm{\alpha \mapsto z(\alpha,\psi,\theta) - \alpha}_{C^k} < \varepsilon
    \qquad\textrm{on } \D_r \times S^1 \times U,
  \end{equation*}
  which completes the proof.
\end{proof}

\subsubsection{Asynchronous Equilibrium in the closed Equation }
\label{sec:closed_equation_alpha}

The Watanabe--Strogatz equations with functions $f,g$ given
by~\eqref{eq:WS-functions-fg} depend on $\alpha,\psi$ (and $\theta$ as
constant of motion) through the function $z(\alpha,\psi,\theta)$
in~\eqref{eq:order-param-func}. This makes it difficult to find fixed
points and analyse their stability. Our approach is to first make the
assumption that $z = \alpha$, thereby `closing the equation' to a
simple form. Secondly, after we find a stable fixed point, we use
Lemma~\ref{lem:z-alpha-close} to prove that it persists after
reinserting the original function $z(\alpha,\psi,\theta)$.

Let us therefore start analyzing the
equations~\eqref{eq:mobius_evolution} with $f$ and $g$ given
by~\eqref{eq:WS-functions-fg} but $z$ replaced by $\alpha$. In this
case the system reduces to a skew-product flow; the equation for
$\alpha$ does not depend on $\psi$ anymore and becomes
\begin{equation}\label{eq:dotalpha-closed}
  \dot{\alpha} =
     -\frac{\sigma}{2}\Bigl(e^{-i \delta} + \beta e^{i \delta}\Bigr)\alpha^2
     + i(1-\beta)  \alpha
     + \frac{\sigma}{2} \Bigl(\beta \abs{\alpha}^2 e^{-i \delta} + e^{i \delta}\Bigr).
\end{equation}
with $0 < \delta < \pi/4$. The next proposition characterises the asynchronous branch.

\begin{theorem}\label{thm:hyperbolicity_incoherent_branch}
Let
\begin{equation}\label{eq:stable-async-cond}
  0 < \sigma < \sigma_f := \frac{\beta-1 }{\sqrt{1+2\beta\cos(2\delta)}}.
\end{equation}
Then equation~\eqref{eq:dotalpha-closed} has an exponentially
attracting fixed point $\alpha^I$ given by
\begin{equation}\label{eq:stable-async-alpha}
  \abs{\alpha^I} = \frac{\beta - 1 - \sqrt{(\beta-1)^2 - \sigma^2(1+2\beta \cos(2\delta))}}{\sigma(1 + 2\beta \cos(2\delta))}
  \quad\text{and}\quad \arg(\alpha^I) = -\frac{\pi}{2} + \delta.
\end{equation}
Moreover, $\alpha^I$ is the unique branch of fixed points satisfying $\lim_{\sigma\rightarrow 0} \alpha^I = 0$.
\end{theorem}

\begin{proof}
Representing equation~\eqref{eq:dotalpha-closed} in polar coordinates,
$\alpha = r e^{i\eta}$, we obtain
\begin{subequations}\label{eq:WS-polar}
  \begin{align}
    \dot{r}    &= \frac{\sigma}{2}(1 - r^2) \cos(\eta-\delta), \label{eq:dotr} \\
    \dot{\eta} &= 1-\beta  - \sigma\beta r \sin(\eta+\delta)
                  -\frac{\sigma}{2}\frac{1+r^2}{r}\sin(\eta-\delta). \label{eq:doteta}
  \end{align}
\end{subequations}
Because of the product structure of equation~\eqref{eq:dotr}, we readily
see that fixed points correspond to either $r=1$ or $\cos(\eta-\delta) = 0$.
Since we are looking for the branch $\alpha^I$ of fixed points that converge
to zero as $\sigma \to 0$, we discard the solution $r=1$.
Hence we substitute $\cos(\eta-\delta) = 0$ and $\sin(\eta-\delta) = \pm 1$
in equation~\eqref{eq:doteta} and obtain the condition
\begin{equation*}
  \pm(\beta-1)  + \sigma\beta r \cos(2\delta) +\frac{\sigma}{2}\frac{1+r^2}{r} = 0.
\end{equation*}
Solving for $r$, we note that $\sin(\eta-\delta) = 1$ leads to negative
solutions, so we find that $\eta = -\pi/2 + \delta$ and
\begin{equation}\label{eq:stable-sync-r}
  r_{\pm} = \frac{(\beta-1)  \pm \sqrt{(\beta-1)^2 - \sigma^2(1+2\beta \cos(2\delta))}}{\sigma(1+2\beta \cos(2\delta))}.
\end{equation}
These fixed points exist for
\begin{equation*}
  \sigma < \sigma_f = \frac{(\beta-1) }{\sqrt{1+2\beta \cos(2\delta)}},
\end{equation*}
where the square root is well-defined since $\delta \in (0,\frac{\pi}{4})$.
The branch $r_-$ uniquely satisfies the condition that it converges to zero
as $\sigma \to 0$.

Next we analyze its stability. Evaluating the Jacobian matrix of the
system~\eqref{eq:WS-polar} at the fixed point we obtain
\begin{equation}\label{eq:jacobian_equation}
  J(r_-,-\frac{\pi}{2}+\delta) =
  \left(
    \begin{array}{cc}
      0 & \frac{\sigma}{2}(1-r_-^2) \\
      a & -\beta \sigma r_- \sin(2\delta)
    \end{array}
  \right)
  \qquad\text{with}\quad
  a = \frac{\sigma}{2}\left[1+2\beta \cos(2\delta) - \frac{1}{r_{-}^2} \right].
\end{equation}
Considering~\eqref{eq:stable-sync-r} and interpreting
$(\beta-1) $ and $\sigma\sqrt{1 + 2\beta\cos(2\delta)}$ as the
sides of a triangle, we apply a triangle inequality estimate and
conclude that
\begin{equation*}
  1+2\beta \cos(2\delta) < \frac{1}{r_-^2}.
\end{equation*}
Therefore, it follows that for $\delta \in (0,\frac{\pi}{4})$ and
$\sigma>0$ we have $a < 0$ and
\begin{equation*}
  \operatorname{tr}(J) < 0 \quad\text{and}\quad \det(J) > 0,
\end{equation*}
hence both eigenvalues of $J$ have negative real part, which completes
the proof.
\end{proof}

\subsubsection{Below backward critical coupling}
\label{sec:sync_diagram_below_backward_critical_coupling}

We are interested in the large $\beta$ regime. In this case,
$\sigma_c^b \rightarrow1$ as $\beta \rightarrow \infty$. The Jacobian
eigenvalues at the incoherence state $\alpha^{I}$ are complex and their
imaginary part is proportional to $\beta$ whereas their real parts
converge to a constant depending on $\delta$ and $\sigma$. As a
consequence, near the fixed point solutions spiral towards the
fixed point with high frequency. We can also construct a Lyapunov
function near $\alpha^I$ but its domain is rather small. This shows
that the nonlinear picture for $|\alpha| \approx 0$ is slightly more
intricate. To better handle the dynamics of $\alpha$ we average
out the fast oscillatory behaviour and analyze the averaged system.

\begin{theorem}\label{thm:Asymptrho}
Let $z$ be the order parameter of the system~\eqref{eq:order-param-func}
corresponding to initial conditions in an open neighborhood $U$ of
 $\Theta \subset \mathbb{T}^{N}$. Given $k\in \mathbb{N}$,
 $\delta_0>0$, $0< \delta< \pi/4$, $0< \sigma<1$ and  $\varepsilon>0$
 small enough, there is $N_0 \in \mathbb{N}$, $\beta_0>0$  such
 that for all $N>N_0$, $\beta> \beta_0$,  and $z(0) \in
 \mathbb{D}_{\delta_0}$  there is $T_0>0$ such that in $\tau := \beta t$
\begin{equation*}
|z(\tau) - \alpha^I | \le \varepsilon, \mbox{~ for all ~ } \tau > T_0.
\end{equation*}
\end{theorem}

The proof of this theorem is presented below. Roughly speaking,
it means that if we follow the synchronized manifold $M_C$
decreasing slowly the intra-coupling strength once it looses stability
the order parameter will drop and stay near zero. Thus, only the
asynchronous branch is stable for $0<\sigma<1$ and large $\beta$.
We will obtain this result as a consequence of three Lemmas. To
this end we first need some preliminary setting.

The planar system~\eqref{eq:WS-polar} in polar coordinates can
be rewritten in the form
\begin{equation*}
\dot{r} = F(r, \eta) \quad  \mathrm{and}  \quad \dot{\eta} = G(r, \eta)
\end{equation*}
where
\begin{subequations}\label{eq:WS-polar_average}
  \begin{align}
    F(r, \eta)    &:=  \frac{\sigma}{2}(1 - r^2) \cos(\eta-\delta), \\
    G(r, \eta)    &:= -\beta [1 + \sigma r \sin(\eta + \delta)]  + \left(1
                  -\frac{\sigma}{2}\frac{1+ r^2}{r}\sin(\eta -\delta) \right).
  \end{align}
\end{subequations}
 We can parametrize time using the parameter $\beta$ and introduce an averaged system on the new time $\tau := \beta t$. First, we need some auxiliary results.

\begin{lemma}\label{lemma:Fbeh}
Fix $\delta_0 \in (0, 1)$ and consider $R : [\delta_0 -1, 1 - \delta_0] \rightarrow \mathbb{R}$ as
\begin{equation*}
  R(x) = \frac{P(x)}{Q(x)},
\end{equation*}
where
\begin{equation*}
  P(x) = - \frac{1}{\pi} \int_0^{2\pi} \frac{\sin \vartheta}{1+ x \sin\vartheta} d\vartheta \mbox{~ ~ ~ }  and  \mbox{~ ~ ~ } Q(x) = \frac{1}{2\pi} \int_0^{2 \pi} \frac{d \vartheta}{1 + x \sin(\vartheta)}.
\end{equation*}
Then $R$ is a smooth function satisfying
\begin{equation*}
  R(0) = 0, \mbox{~~~} R(x) > 0  \mbox{~ for ~} x\in(0,1-\delta_0],  \mbox{ ~ and ~} \frac{dR(0)}{dx} =  1.
\end{equation*}
\end{lemma}
\begin{proof} First we calculate $P(x)$. We Taylor expand the function $1/(1+x \sin \vartheta)$ to obtain
\begin{equation*}
  P(x) = - \frac 1 \pi \int_{0}^{2 \pi}\sin {\vartheta} \sum_{n=0}^{\infty} (-1)^n (x \sin {\vartheta})^n  d\vartheta
\end{equation*}
notice that $\int_0^{2\pi} \sin^{2 k +1} \vartheta d\vartheta = 0$ as $\sin^{2k +1} \vartheta$ is an odd function in $[0,2 \pi]$ and thus
\begin{equation*}
  - \frac{1}{\pi} \int_{0}^{2 \pi} \frac{\sin({\vartheta})}{1 + x \sin({\vartheta})}d\vartheta =
\sum_{k=1}^{\infty} a_{2k} x^{2k - 1}
\end{equation*}
where
\begin{equation*}
  a_{2k} = \frac{1}{\pi}\int_0^{2\pi } \sin^{2k} \vartheta d\vartheta >0, \,\,\,\,\, k\ge 1
\end{equation*}
so $P$ is a monotone function. Moreover,
\begin{equation*}
  a_2 = 1, \mbox{~ ~ ~ } a_4 = \frac{3}{4}, \mbox{~ ~ ~ }  and  \mbox{~ ~ ~ } a_6 = \frac{5}{8}.
\end{equation*}
Repeating the same procedure
\begin{align*}
    Q(x) &= \frac{1}{2 \pi} \int_{0}^{2 \pi} \sum_{n=0}^{\infty} (-1)^n (x \sin {\vartheta})^n d\vartheta \\
        &= 1 + \sum_{k = 1}^{\infty} \frac{a_{2k}}{2} x^{2 k} \\
        &= 1 + \frac{x}{2} P(x).
\end{align*}
Therefore, $R(x)$ is positive, and the claim follows.
\end{proof}

This Lemma will play a role in the averaging approximation for the dynamics of $\alpha$

\begin{lemma}\label{lemma:averaged_system}
   Fix $0<\delta_0 < 1$  and consider $\mathbb{D}_{\delta_0} : = \{ z \in \mathbb{C} : |z| < 1 - \delta_0 \}$. For any $0 < \sigma < 1$, $0<\delta<\pi/4$, $\alpha_0 \in \mathbb{D}_{\delta_0}$ there exist constants $c > 0$, and $\beta_0 > 0$ such that for each $\beta > \beta_0$  the solutions $\alpha(t)$ of Equation~\eqref{eq:WS-polar_average} in $\tau := \beta t$ with initial condition $\alpha_0 \in \mathbb{D}_{\delta_0}$ satisfy
    \begin{equation}\label{eq:error_estimate_average}
        \|r(\tau) - \rho(\tau)\| < \frac{c}{\beta} \qquad \mathrm{for}~ \tau \geq 0,
    \end{equation}
    where $\rho(\tau)$ is the trajectory of the averaged system
    \begin{equation}\label{eq:truncated_averaged_system}
       \rho^{\prime} = \frac{1}{\beta} F^1(\rho)
       \hspace{0.5cm}
    \end{equation}
  with initial condition $\rho(\eta_0) = r_0$,   where $F^1: [\delta_0 -1,1 - \delta_0] \rightarrow \mathbb{R}$ is smooth and satisfying
  \begin{equation*}
    F^1(0) = 0, \mbox{~ ~  ~ ~ } F^1(\rho) < 0 \mbox{ ~ for ~ } \rho\in(0,1-\delta_0]  \mbox{~ ~ and ~ ~ } \frac{d F^1(0)}{d\rho} =  - \frac{1}{4} \sigma \sin 2\delta.
  \end{equation*}
\end{lemma}

\begin{proof}
    Performing a change of time scale $t \mapsto \tau = \beta t$,  we obtain:
    \begin{align}\label{eq:slow_fast_WS-polar}
     \begin{matrix}
        r^{\prime} = \varepsilon~F(r, \eta) \quad \mathrm{and} \quad
        \eta^{\prime} =
        -1 - \sigma r \sin(\eta + \delta) + \varepsilon ~G_1(r, \eta),
 \end{matrix}
    \end{align}
where $\prime$ denotes differentiation with respect to $\tau$ and $\varepsilon = 1/ \beta$. The transformation induces a slow-fast structure where the small parameter is $\varepsilon$. Notice that  $r=1$  is a fixed point of $r$ and thus solutions starting with $r<1$ cannot have $r>1$ and neither have $r=1$ in finite time. Thus, $\eta^{\prime}$ is bounded away from zero and we can use $\eta$ as our new time. To apply the averaging principle, we define the averaged vector field as
    \begin{equation}\label{eq:averaged_system}
        F^1(r) := \frac{\frac{1}{2 \pi} \int_0^{2 \pi } H(r, \vartheta, 0)~d\vartheta}{\frac{1}{2 \pi} \int_0^{2 \pi } \frac{d\vartheta }{- 1 - \sigma r \sin{(\vartheta + \delta)}} },
    \end{equation}
    where
    \begin{equation}\label{eq:r_eta_vector_field}
       H(r, \eta, \varepsilon) =  \frac{\frac{\sigma}{2}(1 - r^2) \cos(\eta-\delta)}{- [1 + \sigma r \sin(\eta + \delta)]  + \varepsilon \left(1
                  -\frac{\sigma}{2}\frac{1+ r^2}{r}\sin(\eta -\delta) \right)} .
    \end{equation}
    The denominator of Equation~\eqref{eq:averaged_system} resembles the $Q(x)$ function of Lemma~\ref{lemma:Fbeh}. So, replacing Equation~\eqref{eq:r_eta_vector_field} into Equation~\eqref{eq:averaged_system}, we calculate the numerator:
    \begin{align*}
        \int_0^{2 \pi} \frac{ \cos(\eta-\delta)}{1 + \sigma r \sin(\eta + \delta)} ~d\eta = \cos({2 \delta}) \int_{\delta}^{2 \pi + \delta} \frac{\cos(\vartheta)}{1 + \sigma r \sin({\vartheta})}d\vartheta + \sin({2\delta}) \int_{\delta}^{2 \pi + \delta} \frac{\sin({\vartheta})}{1 + \sigma r \sin({\vartheta})}d\vartheta,
    \end{align*}
where we used the change of variables $\vartheta = \eta + \delta$ and the trigonometric relation for the cosine. Since both integrals are along a full cycle of $\vartheta$, the first integrals is zero. Then
 \begin{align*}
     -   \frac 1 \pi \int_0^{2 \pi} \frac{ \cos(\eta-\delta)}{1 + \sigma r \sin(\eta + \delta)} ~d\eta =  \sin({2\delta}) P(\sigma r),
    \end{align*}
From this observation we obtain the  vector field for $\rho$. Thus we obtain
\begin{equation*}
F^1(\rho) = - \frac{\sin 2\delta}{4} \sigma (1-\rho^2) R(\sigma \rho).
\end{equation*}
    From Lemma~\ref{lemma:Fbeh} it follows that since $R(\sigma \rho)>0$ for $\rho>0$ we obtain $F^1(\rho)<0$.

To conclude the stability, notice that since $F^1<0$ for $\rho<1$ we have that solutions converge to the origin. Solutions will enter a ball of radius $\delta_1$ sufficiently small.  We show that this convergence is exponential as $\frac{dF^{1}(0)}{d\rho} <0$. The result follows by the principle of linearized stability: given $\delta_1>0$ small enough, all solutions with $\rho < \delta_1$ converge to the origin exponentially fast. Next applying averaging Theorem~\ref{thm:averaging_attractor} our claim follows.
\end{proof}

\begin{lemma}\label{ConvAlpha}
Fix $0<\delta_0 < 1$. Then, for any $0 < \sigma < 1$, $0<\delta<\pi/4$, $\alpha_0 \in \mathbb{D}_{\delta_0}$,  there exists $\beta_0 > 0$ such that for each $\beta > \beta_0$  the solutions $\alpha(t)$ of Equation~\eqref{eq:WS-polar_average} with initial condition $\alpha_0 \in \mathbb{D}_{\delta_0}$ converge to $\alpha^I$.
\end{lemma}

\begin{proof}
Given $\varepsilon_0>0$ small enough, Lemma~\ref{lemma:averaged_system}
shows that there is constant $c$ and $\beta_0$ such that for every $\beta> \beta_0$ we have
$\| \alpha(\tau) \| = \|r(\tau) \| \le \varepsilon_0$, since $0$ is a hyperbolic fixed point for $\rho$.

Equation~\eqref{eq:WS-polar_average} has a stable fixed point $\alpha^I$ satisfying $|\alpha^I| \rightarrow 0$ as $\beta \rightarrow \infty$ as can be observed in Equation~\eqref{eq:stable-sync-r}. Moreover, when $\beta\rightarrow \infty$ the real part of the Jacobian eigenvalues of $\alpha^I$ tends to a constant independent of $\beta$. Thus, by the principle of linearized stability there is $\delta_I>0$ independent of $\beta$ such that any initial condition in the open ball $B(\alpha^I,\delta_I)$ converges to $\alpha^I$. Now we choose $\varepsilon_0 < \delta_I$. This is guaranteed as $\delta_I$ is independent of $\beta$.
Thus, when $\beta$ is large enough, solutions $\alpha$ starting in $\mathbb{D}_{\delta_0}$ enter $B(\alpha^I,\delta_I)$ and thus must converge to $\alpha^I$.
\end{proof}

\begin{proof}[Proof of Theorem~\ref{thm:Asymptrho}]
By Lemma~\ref{lem:z-alpha-close} there is $N_0$ such that for all $N > N_0$
\begin{equation*}
|  z(\tau) - \alpha (\tau) | \le \varepsilon/2.
\end{equation*}
Next, by the averaging principle there exists $c$ and $\beta_0$ such that for all $\beta > \beta_0$ we find constants $M(\beta)$ and $\mu(\beta)$ such that
\begin{align*}
| \alpha(\tau) - \alpha^I| &\leq |r(\tau) - \rho(\tau)| + 2 |\rho(\tau)| + |\rho(\tau) - r_{-}| \\
    &  \leq \frac{c}{\beta} + 2 M ~|\rho_0| ~e^{-\mu \tau} + M ~|\rho_0 - r_{-}|~e^{-\mu \tau}.
\end{align*}
Consequently, after a fixed finite time $T_0$ we can make it below $\varepsilon/2$. Applying the triangle inequality for $t > T_0$ we obtain
\begin{equation*}
| z(t) - \alpha^I| \le | z(t) - \alpha(t) | + | \alpha(t) - \alpha^I|
\end{equation*}
and the claim follows.
\end{proof}

The final step is to construct the invariant manifold relative to the asynchronous dynamics. To this end we use a persistence argument over the equation of $\alpha$ and $\psi$.

\begin{proposition}\label{prop:async_state_NAIM}
  Let $\sigma$ strictly satisfy
  inequality~\eqref{eq:stable-async-cond}. Given $\varepsilon >0$ there
  exist $N_0$ such that for each $N \ge N_0$ there exists a normally
  attracting invariant manifold with
  invariant boundary, $M_I$, such that the order parameter $z$ on
  $M_I$ satisfies
  \begin{equation*}
    \abs{ z - \alpha^I } \le \varepsilon,
  \end{equation*}
  where $\alpha^I$ is given by~\eqref{eq:stable-async-alpha}.
\end{proposition}

\begin{proof}
To construct the incoherent invariant
manifold, we change to coordinates $\alpha,\psi$ using
the Watanabe--Strogatz approach~\ref{sec:ws_theory} while viewing the
$\theta_i$'s as fixed parameters.
First consider the modified system~\eqref{eq:mobius_evolution} with
$\alpha$ as argument of the functions $f,g$. Then the equation for
$\alpha$ decouples and we find a stable fixed point
$\alpha^I \in \C$, with $\abs{\alpha^I} < 1$ and explicitly given by~\eqref{eq:stable-async-alpha}.

Since the dynamics for $\psi$ is neutral, it follows that
$M_\text{WS} = \{ \alpha = \alpha^I, \; \psi \in S^1 \} \subset \D \times S^1$
is a NAIM. Lemma~\ref{lem:z-alpha-close}
shows that given a $k \ge 1$ and a ball $B_r(\alpha^I) \subset \D$, there
exists an $N_0$ and for all $N \ge N_0$ a neighborhood
$U \subset \T^N$ of the uniformly distributed states $\Theta$, such
that the substitution of $\alpha$ by $z(\alpha,\psi,\theta)$ is a
$C^k$-small perturbation of Eq.~\eqref{eq:mobius_evolution}.

We choose $k = 1$ and by persistence of NAIMs,
see~\cite{fenichel1971,HPS1977}, the original system with the true
dynamics for $z$ also has a NAIM
\begin{equation}\label{eq:invar-incoherent}
  \tilde{M}_\text{WS}(\theta) = \{ \alpha = \alpha^I + g^I(\psi,\theta), \;
  \psi \in S^1 \}
  \subset \D \times S^1
\end{equation}
$C^1$-close to $M_\text{WS}$. Thus, since $U$ is precompact, we have
uniformly for all $\theta\in U$ that $\norm{g^I}_{C^1} \le \varepsilon$ when
$N_0$ is sufficiently large, and thus on $\tilde{M}_\text{WS}(\theta)$
also $\abs{z - \alpha^I} \le \varepsilon$ holds. Now we consider
$\theta \in \T^N$ as dynamical variables again; that is, we lift the
system according to~\eqref{eq:diagram-WS-coords} to coordinates
$(\alpha,\psi,\theta)$. Since the $\theta_i$'s have no dynamics, also
\begin{equation}\label{eq:invar-incoherent-full}
M_I = \bigcup_{\theta \in U} \tilde{M}_\text{WS}(\theta) \times
\{\theta\} \subset \D \times S^1 \times \T^N
\end{equation}
is a NAIM, but with
boundary\footnote{%
  Here we mean the topological boundary, that is,
  $\partial M_I = \overline{M}_I \setminus M_I$.%
}. Since $\theta$ has no associated dynamics, the boundary is
invariant.

\end{proof}

\section{Statement and Proof of the Main Theorem}
\label{sec:main-proof}

In this section, we state and prove the main result of this paper~\ref{thm:chimera-state}. First we must recapitulate the results we obtained in Section~\ref{sec:mathematical_analysis_prelimanaries_and_results}.

Consider two star networks with a large number of leaves $N$. The first star we shall assume to be in a synchronized, coherent state (labeled $+$ and with coordinates labeled with $^+$) while the second network we assume to be in an incoherent state (labeled $-$ and with coordinates labeled with $^-$). Both the coherent and incoherent states correspond to normally attracting invariant submanifolds of each star network:

\begin{itemize}
    \item The coherent invariant manifold is given in $\phi_i^+,\varphi_0^+$
coordinates with $i = 1, \dots, N$ by
    \begin{equation*}
        M_+ = \{ \phi_1^+ = \ldots = \phi_N^+ = \phi^C, \varphi^+_0 \in S^1 \}
    \end{equation*}

        where $\phi^C$ is fixed and given by~\eqref{eq:stable-sync-phi}.
    \item The incoherent invariant manifold is given in $\varphi^-_0, \alpha,\psi,\theta \in Q^{-} := S^1 \times \mathbb{D} \times S^1 \times \mathbb{T}^N$ coordinates by

    \begin{equation*}
        M_- = \{ \alpha = \alpha^I + g^I(\psi,\theta),\psi \in S^1, \theta \in U \subset \T^N, \varphi^-_0 \in S^1 \}
    \end{equation*}

    with $\alpha^I$ given by~\eqref{eq:stable-async-alpha} and $g^I$ a $C^1$-small function for $N$ sufficiently large, following Proposition~\ref{prop:async_state_NAIM}.
\end{itemize}

    Note that the product system has $M = M_+ \times M_-$ as normally attracting invariant manifold. Moreover, the order parameters $z^{+}_{M}:\mathbb{T}^{N + 1} \to \mathbb{C}$ and $z^{-}_{M}: Q^{-} \to \mathbb{C}$ are understood as the restriction on $M$ and this is extended when we couple the two stars with intra-coupling $\varepsilon$ guaranteeing the existence of chimera states as presented in Theorem~\ref{thm-NAIM}.


\begin{theorem}[Normally Attracting Invariant Manifold for weakly coupled stars]
\label{thm-NAIM}
  Consider the coupled star network system~\eqref{eq:coupled_stars}.  For any $h \in C^1$ and $\eta>0$, there exist $\varepsilon_0>0$, $N_0>0$, $\beta_0>0$ and an open set $\mathcal{I} \subset \mathbb{R}$ such that for any $N > N_0$, $\beta > \beta_0$ and $0 < \varepsilon \le \varepsilon_0$, $\sigma \in \mathcal{I}$ the system has a normally
    attracting invariant manifold with boundary $M_{\varepsilon}$ that is
    $\mathcal{O}(\varepsilon)$-close to $M = M_+ \times M_- \subset \mathbb{T}^{N + 1} \times Q^{-}$ such that
\begin{equation*}
   r^{+}_{M_{\varepsilon}}(\varphi_0^{+}, \phi_1^{+}, \dots, \phi_N^{+})   >  1 - \eta
\end{equation*}
and
\begin{equation*}
  r^{-}_{M_{\varepsilon}}(\alpha, \psi, \theta, \varphi_0^{-})  < \eta.
\end{equation*}
\end{theorem}
The choice of coupling function $h$ between each star affects the order of time these chimera states exist. So, we split the results in Theorem~\ref{thm:chimera-state} and Theorem~\ref{thm:asymptotically_stable_chimera-state}.

\begin{theorem}[At least metastable Chimera States]\label{thm:chimera-state}
\leavevmode 
\begin{enumerate}
    \item
    \label{enum:thm-eps-time}\emph{[General Coupling]}
    If $h \in C^1$ in Theorem~\ref{thm-NAIM} is arbitrary, there exists an open set $U \subset \mathbb{T}^{2(N + 1)}$
    of initial conditions $\theta \in U$
    on $M_{\varepsilon}$ such that the coupled dynamics stays in $M_{\varepsilon}$ for a time of
    order $1/\varepsilon$, uniformly for $\theta \in U$. In other words, chimera states surely exist for a time of order $1/\varepsilon$
  \item\label{enum:thm-eps2-time}\emph{[Kuramoto--Sakaguchi Coupling]} In particular, if the coupling $h \in C^{\infty}$ is of the form
    \begin{equation*}
      h(\phi_1,\phi_2) = c_1 \sin(\phi_1-\phi_2 + \vartheta).
    \end{equation*}
    Then there exists an open set $U \subset \mathbb{T}^{2(N + 1)}$
    of initial conditions $\theta \in U$
    on $M_{\varepsilon}$ such that the coupled dynamics stays in $M_{\varepsilon}$ for a time of
    order $1/\varepsilon^2$, uniformly for $\theta \in U$.
\end{enumerate}
\end{theorem}

If extra information on the coupling is given or the coupling structure is not fully symmetric, then the chimeras can be asymptotically stable. This motivates us to highlight and discuss how general is our results:
\begin{remark}
 Let $A \in \mathbb{R}^{N + 1\times N + 1}$ be the adjacency matrix corresponding to the intercoupling between both stars. In a general flavor the networks dynamics we analyse is given as
    \begin{equation}\label{eq:coupled_stars_generalization}
        \begin{split}
            \dot{\varphi}_0^+ &= \beta + \beta \frac{\sigma}{N} \sum_{j=1}^{N} \sin(\varphi_j^+ - \varphi_0^+ + \delta)
                          + \varepsilon \sum_{j=0}^N A_{0j} \, h(\varphi^+_0 - \varphi^-_j),\\
            \dot{\varphi}_i^+ &= 1 + \sigma \sin(\varphi_0^+ - \varphi_i^+ + \delta)
                          + \varepsilon \sum_{j=0}^N A_{ij} \, h(\varphi^+_i - \varphi^-_j), \quad i = 1,\dots, N,
        \end{split}
    \end{equation}
    and likewise for the ``$-$" star. Both previous theorems, Theorem~\ref{thm-NAIM} and Theorem~\ref{thm:chimera-state}-(\ref{enum:thm-eps-time}), are extended to this general coupling structure between stars without any extra condition, since the persistence argument goes through.
\end{remark}
  In particular, for the case of small number of connections between both stars, we state Theorem~\ref{thm:asymptotically_stable_chimera-state}.

 \begin{theorem}[Stable Chimera States]\label{thm:asymptotically_stable_chimera-state}
 Consider the coupled star network system~\eqref{eq:coupled_stars_generalization}.
 For any $h \in C^1$ $\eta>0$ there exists $\varepsilon_0 > 0$, $N_0 > 0$, $\beta_0 >0 $ and an open set $\mathcal{I} \subset \mathbb{R}$ such that for any $\beta > \beta_0$, $N > N_0$ and $0 \le k/N \ll 1$, where only $k$ of the leaves are
    coupled, that is, $A_{ii} = 1$ for $i \in I$ and $A_{ij} = 0$
    otherwise, with $\abs{I} = k$.  Then for any $0 < \varepsilon \le \varepsilon_0$, $\sigma \in \mathcal{I}$ and  initial conditions on
    $M_\varepsilon$ the system stays in a product of coherent and incoherent
    states for all time, even though it may leave $M_{\varepsilon}$.
\end{theorem}
\subsection{NAIM for Coupled Stars}

When we couple the two stars we introduce a coupling of size $\varepsilon$ between the
two systems, with $\varepsilon \le \varepsilon_0$ sufficiently small, then
this perturbation of the product system will again have a normally
attracting invariant manifold $M_\varepsilon$ close to $M$. This manifold
can be described as
\begin{equation}\label{eq:coupled-invar-mfld}
  M_\varepsilon = \Bigl\{\,
                \phi_i^+ = \phi^C + \tilde{g}_i^C(\varphi^\pm_0,\psi,\theta),\;
                \alpha   = \alpha^I + (g^I+\tilde{g}^I)(\varphi^\pm_0,\psi,\theta),\;
                \varphi^\pm_0, \psi \in S^1,\; \theta \in U \subset \T^N
              \,\Bigr\},
\end{equation}
with functions $\norm{\tilde{g}^C}_{C^1},\norm{\tilde{g}^I}_{C^1} \in \mathcal{O}(\varepsilon)$
describing the perturbation of $M_\varepsilon$ away from $M$. Note that
since the combined system is invariant under a global phase shift, the
functions $\tilde{g}^C,\tilde{g}^I$ only depend on $\varphi^\pm_0$
through their phase difference $\varphi^+_0 - \varphi^-_0$.
Finally, since $M_\varepsilon$ is a NAIM and a normally
hyperbolic invariant manifold, it follows that $M_\varepsilon$ has an
invariantly foliated stable manifold $W^s(M_\varepsilon)$,
see~\cite[Thm.~4.1]{HPS1977}. This implies that the orbits of points
in a leaf $W^s(m)$ shadow and exponentially in time converge to the
orbit of $m \in M_\varepsilon$.

\begin{figure}[htb]
  \centering
  \input{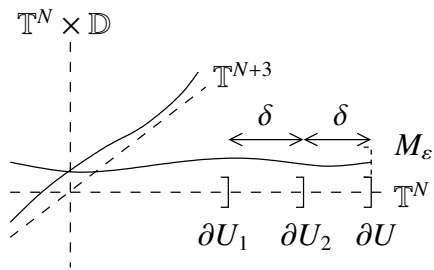}
  \caption{Flow along the invariant manifold $M_\varepsilon$.}
  \label{fig:boundary-M}
\end{figure}

Although $M_\varepsilon$ has a boundary, both the persistence result and
existence of the invariant stable foliation hold true on $M$ with an
arbitrarily small neighborhood of its boundary removed. That is, we
modify the perturbed dynamics in a vertical neighborhood over
$\partial U$, such that the modified coupled dynamics leaves a
vertical section $S= \T^N \times \D \times \T^{N+3} \times \partial U$
over $\partial M$ invariant. Hence persistence of $M$ to $M_\varepsilon$
is well-defined, where $M_\varepsilon$ again has an invariant boundary
with $\partial M_\varepsilon \subset S$, see~Figure~\ref{fig:boundary-M}.
This follows from Theorem~4.8 in~\cite[Sec.~4.2]{eldering2013}. Since
$M_\varepsilon$ is invariant, it again has an invariant stable foliation.
As long as solutions of the coupled system on $M_\varepsilon$ stay away
from $\partial M_\varepsilon$, they satisfy the unmodified dynamics and
the conclusions drawn above remain true. This proves Theorem~\ref{thm-NAIM}.

\subsection{General Coupling}

For point~\textit{\ref{enum:thm-eps-time}} we have to study how the
coupling may perturb the trivial dynamics $\dot{\theta} = 0$ in the
incoherent system, and drive the system across the boundary of $U$,
where normal hyperbolicity subsequentially may be lost. In view
of~\eqref{eq:coupled-invar-mfld}, we can express the invariant
manifold $M_\varepsilon$ for the coupled system as a graph, where the
coordinates $\phi_i^+$ and $\alpha$ depend on the other coordinates
$\varphi_0^+,\varphi_0^-,\psi,\theta \in \T^{N+3} \times U$, see
Figure~\ref{fig:boundary-M}. Now the boundary of $M_\varepsilon$ sits
over $\T^{N+3} \times \partial U$. Solution curves that
leave $M_\varepsilon$ do so by $\theta$ leaving $U$ through its boundary
$\partial U$.

In the uncoupled setup there was no dynamics along $\theta \in U \subset \T^N$,
that is, $\dot{\theta} = 0$, thus, with coupling we have
$\dot{\theta} \in \mathcal{O}(\varepsilon)$. We choose open
neighborhoods $U_1, U_2$ such that
\begin{equation*}
  \Theta \subset U_1 \subset U_2 \subset U \subset \T^N
  \qquad\text{with } B(U_1,\delta) \subset U_2
        \text{ and } B(U_2,\delta) \subset U
\end{equation*}
for some fixed $\delta > 0$, see also Figure~\ref{fig:boundary-M}.
We may also assume that the modifications we made for Theorem~\ref{thm-NAIM} are contained outside $U_2$.
Solution curves with initial conditions in $M_\varepsilon$ and
$\theta \in U_1$ will take
$\mathcal{O}\bigl(\frac{\delta}{\varepsilon}\bigr)$ time to cross the
gap $U_2 \setminus U_1$. As the vector field is unmodified in this
part of phase space, these conclusions are true for the unmodified
system. This proves point~\textit{\ref{enum:thm-eps-time}}.

\subsection{Full one-to-one sinusoidal coupling}

We continue to prove point~\textit{\ref{enum:thm-eps2-time}}.
Consider Equation~\eqref{eq:coupled_stars} and assume that the coupling is
sinusoidal, that is, $h = c_1 \sin {}+ c_2 \cos$. Thus the equations
for $\phi^-$ become
\begin{equation}\label{eq:async-sin-coupled}
  \dot{\phi}^-_i = (1-\beta)  - \sigma \sin(\phi^-_i - \delta) -\beta \frac{\sigma}{N} \sum_{j=1}^{N} \sin(\phi_j^- + \delta)
                  +\varepsilon\,h(\phi^+_i - \phi^-_i + \varphi^+_0 - \varphi^-_0).
\end{equation}
Since we know already that the Watanabe--Strogatz approach can be applied to
the first terms comprising the uncoupled system, we focus on the coupling
term $h(\phi^+_i - \phi^-_i + \varphi^+_0 - \varphi^-_0)$. We denote
the associated vector field by $X \in \mathfrak{X}(\T^N)$, but note
that $X$ also depends on $\phi^+, \varphi^\pm_0$ as external
variables. We can write this as
\begin{equation*}
  h(\phi^+_i - \phi^-_i + \varphi^+_0 - \varphi^-_0)
  = h_1(-\phi^-_i) h_2(\phi^+_i + \varphi^+_0 - \varphi^-_0),
\end{equation*}
where $h_1,h_2$ again denote sinusoidal functions. We cannot
 apply the Watanabe--Strogatz lift~\eqref{eq:WS-lift} to
$X$, since $h_2(\phi^+_i + \varphi^+_0 - \varphi^-_0)$ is not
identical for all $i$. However, when we restrict to the invariant
manifold $M_\varepsilon$, then
$\phi^+_i = \phi^C+\tilde{g}^C_i(\varphi^+_0-\varphi^-_0,\psi,\theta)$
with $\tilde{g}^C_i \in \mathcal{O}(\varepsilon)$. Hence we find
\begin{equation*}
  h(\phi^+_i - \phi^-_i + \varphi^+_0 - \varphi^-_0)
  = h_1(-\phi^-_i) h_2(\phi^C + \varphi^+_0 - \varphi^-_0) + \mathcal{O}(\varepsilon).
\end{equation*}
Now we can lift the vector field $X'$ associated to
$h_1(\phi^-_i) h_2(\phi^C + \varphi^+_0 - \varphi^-_0)$ using the
Watanabe--Strogatz lift~\eqref{eq:WS-lift}. This implies that
$\hat{X}'$ leaves $\theta$ constant, and since the remaining part of
the coupling vector field is of order $\varepsilon^2$, it follows that
$\dot{\theta} \in \mathcal{O}(\varepsilon^2)$. This proves
claim~\ref{enum:thm-eps2-time} of Theorem~\ref{thm:chimera-state}.

\subsection{Sparse one-to-one coupling}

Finally, we prove Theorem~\ref{thm:asymptotically_stable_chimera-state}. We now
assume that the coupling matrices are of the form $A_{ii} = 1$
when $i \in I \subset \N$ and zero otherwise, where
$\abs{I} = k \ll N$. The coupling function $h \in C^1$ is arbitrary;
since $h$ is defined on a compact set, $\norm{h}_{C^1}$ is bounded.

When only a small fraction of the leaves of the two stars are coupled,
we can improve the previous argument to obtain that the chimera state
is fully stable for all time. To this end, we apply the two different
lifts of the vector field for the incoherent star. We apply the
Watanabe--Strogatz lift~\eqref{eq:WS-lift} to the original uncoupled
vector field and the lift~\eqref{eq:lifted-vectorfield} to the
coupling terms.

Without coupling term, the dynamics of the unsynchronized star are
given in Watanabe--Strogatz lifted coordinates
by~\eqref{eq:mobius_evolution}, \eqref{eq:order-param-func} and
$\dot{\theta} = 0$. The extra coupling term is given in $\phi$
coordinates by
\begin{equation}\label{eq:node-coupling}
  \dot{\phi}_i^\pm = \varepsilon \, h(\phi_i^\mp - \phi_i^\pm \pm \Gamma)
  \qquad\text{for } i \in I,
\end{equation}
where $\Gamma = \varphi_0^+ - \varphi_0^-$ is the difference of the
hub phases. Let us denote by $\varepsilon\,Y \in \mathfrak{X}(\T^N)$
the vector field describing the evolution of the $\phi_i^-$
variables according to~\eqref{eq:node-coupling}. We shall
use~\eqref{eq:lifted-vectorfield} to find the lift
$\hat{Y} \in \mathfrak{X}(\D \times S^1 \times \T^N)$ of $Y$ that
gives equivalent dynamics in Watanabe--Strogatz coordinates. First we
calculate $M_{\alpha,\psi}'(\theta) = \frac{\partial \phi}{\partial
\theta}(\alpha,\psi,\theta)$ with $\phi = \pi(\alpha,\psi,\theta)$.
From $e^{i \phi} = M_{\alpha,\psi}\bigl(w\bigr)$ with
$w = e^{i \theta}$ we obtain
\begin{eqnarray*}
\hphantom{\Longleftrightarrow}\quad\;
  i e^{i \phi} \frac{\partial \phi}{\partial \theta}
  &= \frac{\partial}{\partial w}\Biggl[
       \frac{\alpha + e^{i\psi}w}{1 + \bar{\alpha}e^{i\psi}w} \Biggr] \, i e^{i \theta} \\
\Longleftrightarrow\quad
  \frac{\alpha + e^{i\psi}w}{1 + \bar{\alpha}e^{i\psi}w} \frac{\partial \phi}{\partial \theta}
  &= \Biggl[ \frac{e^{i\psi}}{1 + \bar{\alpha}e^{i\psi}w}
            -\frac{\alpha + e^{i\psi}w}{(1 + \bar{\alpha}e^{i\psi}w)^2}\bar{\alpha}e^{i\psi} \Biggr] w \\
\Longleftrightarrow\quad
  (\alpha + e^{i\psi}w) \frac{\partial \phi}{\partial \theta}
  &= \frac{(1 + \bar{\alpha}e^{i\psi}w)e^{i\psi}w - \alpha\bar{\alpha}e^{i\psi}w - \bar{\alpha}e^{2i\psi}w^2}%
          { 1 + \bar{\alpha}e^{i\psi}w} \\
\Longleftrightarrow\quad
  \frac{\partial \phi}{\partial \theta}(\alpha,\psi,\theta)
  &= \frac{(1 - \alpha\bar{\alpha}) e^{i(\psi+\theta)}}%
          {(1 + \bar{\alpha}e^{i(\psi+\theta)})(\alpha + e^{i(\psi+\theta)})}
   = \frac{1 - \abs{\alpha}^2}{\abs{\alpha + e^{i(\psi+\theta)}}^2}.
\end{eqnarray*}
Thus we find that $\hat{Y}$ is given by
\begin{equation}\label{eq:WS-lifted-coupling}
  \dot{\theta}_i
  = \varepsilon ~\Biggl(\frac{\partial \phi}{\partial \theta}\Biggr)^{-1} \cdot
    h\bigl(\phi_i^- - \phi_i^+ - \Gamma\bigr)
  = \varepsilon ~\frac{\abs{\alpha + e^{i(\psi+\theta_i)}}^2}{1 - \abs{\alpha}^2} \,
    h\bigl(\pi_i(\alpha,\psi,\theta_i) + \varphi_0^- - \varphi_i^+\bigr)
    \quad\text{for } i \in I,
\end{equation}
$\dot{\theta}_i = 0$ for $i \not\in I$, $\dot{\alpha} = 0$ and
$\dot{\psi} = 0$. Although the particular form
of~\eqref{eq:WS-lifted-coupling} should be useful for more detailed
calculations using averaging theory, we did not pursue this. We shall
simply use the fact that the $\theta_i$ with $i \not\in I$ are still
conserved. We use a slightly modified version of
Lemma~\ref{lem:z-alpha-close} to show that if $k \ll N$, then the
result of the Lemma still holds true with the values of $\theta_i$
with $i \in I$ arbitrary, hence there exists an invariant manifold
with invariant boundary for the incoherent star network.

For simplicity, let us assume by relabelling that
$I = \{ N-k+1, \ldots, N \}$. Define the set
\begin{equation}\label{eq:Theta-bar}
  \bar{\Theta} :=
  \bigl\{ (\vartheta,\xi) \in \T^{N-k} \times \T^k \mid
          \exists \theta \in \Theta\colon \forall 1 \le i \le N-k\colon \vartheta_i = \theta_i \bigr\} \subset \T^N.
\end{equation}
Note that $\bar{\Theta}$ is a simple extension of $\Theta$ with the
first $N-k$ entries still uniformly distributed according to
Definition~\ref{def:unif-distrib} and the last $k$ arbitrary. The
coupled dynamics leaves $\bar{\Theta}$ invariant since
$\dot{\vartheta} = 0$, still. We now have the following modified
version of Lemma~\ref{lem:z-alpha-close}; here we fix a $C^1$-norm,
but this is all we needed anyways.

\begin{lemma}\label{lem:z-alpha-close-few-coupled}
  Fix $0 < r < 1$ and consider the disc
  $\D_r = \{ \alpha \in \C \mid \abs{\alpha} < r \}$. For each
  $\varepsilon > 0$ there exists for all $N$ sufficiently large and
  $k \ll N$ an open neighborhood $U = U' \times \T^k$ of\/
  $\bar{\Theta} \subset \T^N$ such that the map
  $\D_r \to \C,\; \alpha \mapsto z(\alpha,\psi,\theta) - \alpha$, with
  $z(\alpha,\psi,\theta)$ given by~\eqref{eq:order-param-func}, is
  smaller than $\varepsilon$ in $C^k$-norm, uniformly for all
  $(\psi,\theta) \in S^1 \times U$.
\end{lemma}

\begin{proof}
  The proof is a straightforward extension of the proof of
  Lemma~\ref{lem:z-alpha-close}. Given
  $\theta' = (\vartheta,\xi) \in \bar{\Theta}$, let
  $\theta \in \Theta$ as in~\eqref{eq:Theta-bar}. For
  $z(\alpha,\psi,\theta')$ we can estimate the difference with
  $z(\alpha,\psi,\theta)$ as
  \begin{equation*}
    z(\alpha,\psi,\theta') - z(\alpha,\psi,\theta)
    = \frac{1}{N} \sum_{j=N-k+1}^N \Biggl(
        \frac{\alpha + e^{i \psi} e^{i \xi_j   }}{1 + \bar{\alpha} e^{i \psi} e^{i \xi_j}}
      - \frac{\alpha + e^{i \psi} e^{i \theta_j}}{1 + \bar{\alpha} e^{i \psi} e^{i \theta_j}}\Biggr).
  \end{equation*}
  Since $\abs{\alpha} < r$ is bounded away from $1$, all terms are
  uniformly bounded in $C^1$-norm, say by $B$. Hence
  $\norm{z(\alpha,\psi,\theta') - z(\alpha,\psi,\theta)}_{C^1} \le
  \frac{k}{N} B$, and thus for $k \ll N$ this yields a sufficiently
  small contribution to the overall $C^1$-norm of
  $\alpha \mapsto z(\alpha,\psi,\theta') - \alpha$ relative to the
  estimate for $z(\alpha,\psi,\theta)$ in the proof
  of~\ref{lem:z-alpha-close}.
\end{proof}

To prove Theorem~\ref{thm:asymptotically_stable_chimera-state}, we now use
Lemma~\ref{lem:z-alpha-close-few-coupled} instead of
Lemma~\ref{lem:z-alpha-close}. Note that $U = U' \times \T^k$ is
invariant under the coupled dynamics since its boundary
$\partial U = \partial U' \times \T^k$ is invariant. Thus, we find
that also the persistent manifold $M_\varepsilon$ has an invariant
boundary.

\subsection{Coupling through Hubs}

Connecting two stars using a general coupling mechanism invalidates the use of the Watanabe--Strogatz approach and the description of the dynamics through the action of the M\"obius group. An important exception to this rule is when only the hubs of the stars are coupled, since, as we will show, we can still apply the WS dimensional reduction method, even if now we have that the phase difference of the hubs acts as an \emph{external field} to the equations describing the order parameters $z$.

Consider two stars with their hubs coupled sinusoidally, with strength $\varepsilon$,
\begin{equation}\label{eq:coupled_stars_hub}
  \begin{split}
    \dot{\varphi}_{k}^{\pm} &= 1 + \sigma \sin(\varphi_{0}^{\pm} - \varphi_{k}^{\pm} + \delta),\\
    \dot{\varphi}_{0}^{\pm} &= \beta + \frac{\sigma \beta}{N} \sum_{j=1}^{N} \sin(\varphi_{j}^{\pm} - \varphi_{H}^{\pm} + \delta) + \varepsilon \sin(\varphi_{0}^{\mp} - \varphi_{0}^{\mp} + \delta),
  \end{split}
\end{equation}
Again we change to coordinates $\phi_{k}^{\pm} = \varphi_{k}^{\pm} - \varphi_{0}^{\pm}$ that measure the relative phases of the leaves with respect to the hubs.
Secondly, we introduce $\Gamma = \varphi_0^+ - \varphi_0^-$, the difference between the hub phases.
With these new variables, the equations become
\begin{align*}
  \dot{\phi}_{k}^{\pm} &= 1 - \beta  - \sigma \sin(\phi_{k}^{\pm} - \delta) - \frac{\beta \sigma}{N} \sum_{j=1}^{N} \sin(\phi_{j}^{\pm} + \delta) + \varepsilon \sin(\mp\Gamma + \delta),\\
  \dot{\Gamma} &= \frac{\sigma \beta}{N} \sum_{j=1}^{N} \left( \sin(\phi_{k}^{+} + \delta) - \sin(\phi_{k}^{-} + \delta) \right) - 2 \varepsilon \cos(\delta) \sin(\Gamma).
\end{align*}
For each star the equations of motion are in the appropriated form and we can apply the WS approach, where the functions $f^{\pm}$ and $g^{\pm}$ are now given by
\begin{align*}
  f^{\pm} &= \frac{-\sigma e^{-i\delta}}{2i},\\
  g^{\pm} &= 1 - \beta  - \beta\sigma \left(\frac{z^\pm e^{i\delta} - \bar{z}^\pm e^{-i\delta}}{2i}\right)
            + \varepsilon \sin(\mp\Gamma + \delta).
\end{align*}
The difference here is that the functions $g^{\pm}$ now depend on the phase difference $\Gamma$ of the hubs, that, in its turn, depends on the relative phases of the leaves of both populations. Nevertheless, the dependence is the same for all leaves and therefore follows the WS approach. The equations for the WS variables $\alpha^{\pm}$ are thus given by
\begin{equation*}
  \dot{\alpha}^{\pm} =
  -\frac{\sigma}{2}(e^{-i\delta} + \beta e^{i\delta})(z^{\pm})^2
  +i(1-\beta)\left( 1 + \varepsilon \frac{\sin(\mp\Gamma + \delta)}{1 - \beta} \right) z^{\pm}
  +\frac{\sigma}{2}\beta e^{-i\delta} |z^{\pm}|^2 + \frac{\sigma}{2}e^{i\delta},
\end{equation*}
where the difference now is that where before would appear only the natural frequency of the leaves $\omega = 1$, now we have a term $1 + \varepsilon \frac{\sin(\mp\Gamma + \delta)}{1 - \beta}$ that comes from the coupling.
The equation for $\Gamma$ can be written as
\begin{equation*}
  \dot{\Gamma} = \sigma \beta \Im( (z^{+} - z^{-})e^{-i\delta}) - 2\varepsilon \cos(\delta) \sin(\Gamma)
\end{equation*}
We can view the effect of the coupling as inducing oscillations of amplitude $\varepsilon/(1-\beta)$ to the natural frequency $\omega$ and therefore we expect the persistence of the fixed points.

\section{Open question and Conclusions}
\label{sec:conclusions}

We presented chimera states born from the coexistence of synchronous and asynchronous dynamics in the mean-field of a single star graph. We showed that the chimera states in coupled star graphs correspond to the existence of an invariant manifold with boundary.  We associated their breaking with dynamics running into the boundary of the invariant manifold and provided lower bounds for their survival. Predicting whether such chimeras break to a complete synchronous dynamics of the network or to asynchronous remains an open question.

The coupling strength between the stars provides the time scale the chimera  exists -- either of order $1/\varepsilon$ (in the general inter-coupling function) or $1/\varepsilon^2$ (for sinusoidal). We performed several numerical experiments to check these predictions. We will discuss this now. To this end, we introduce a parameter $\eta$ which splits the absolute value of the order parameter in Theorem~(\ref{thm-NAIM}) measured on the coherent and the incoherent star. The chimera lifetime $\tau$ is the minimum time it takes to violate the splitting condition from Theorem~(\ref{thm-NAIM}), so:
\begin{equation*}
  |z^- - a^{I}| > \eta \quad \mathrm{or} \quad r^+ < 1 - \eta .
\end{equation*}

\begin{figure}
\centering
\includegraphics[width=0.4\linewidth]{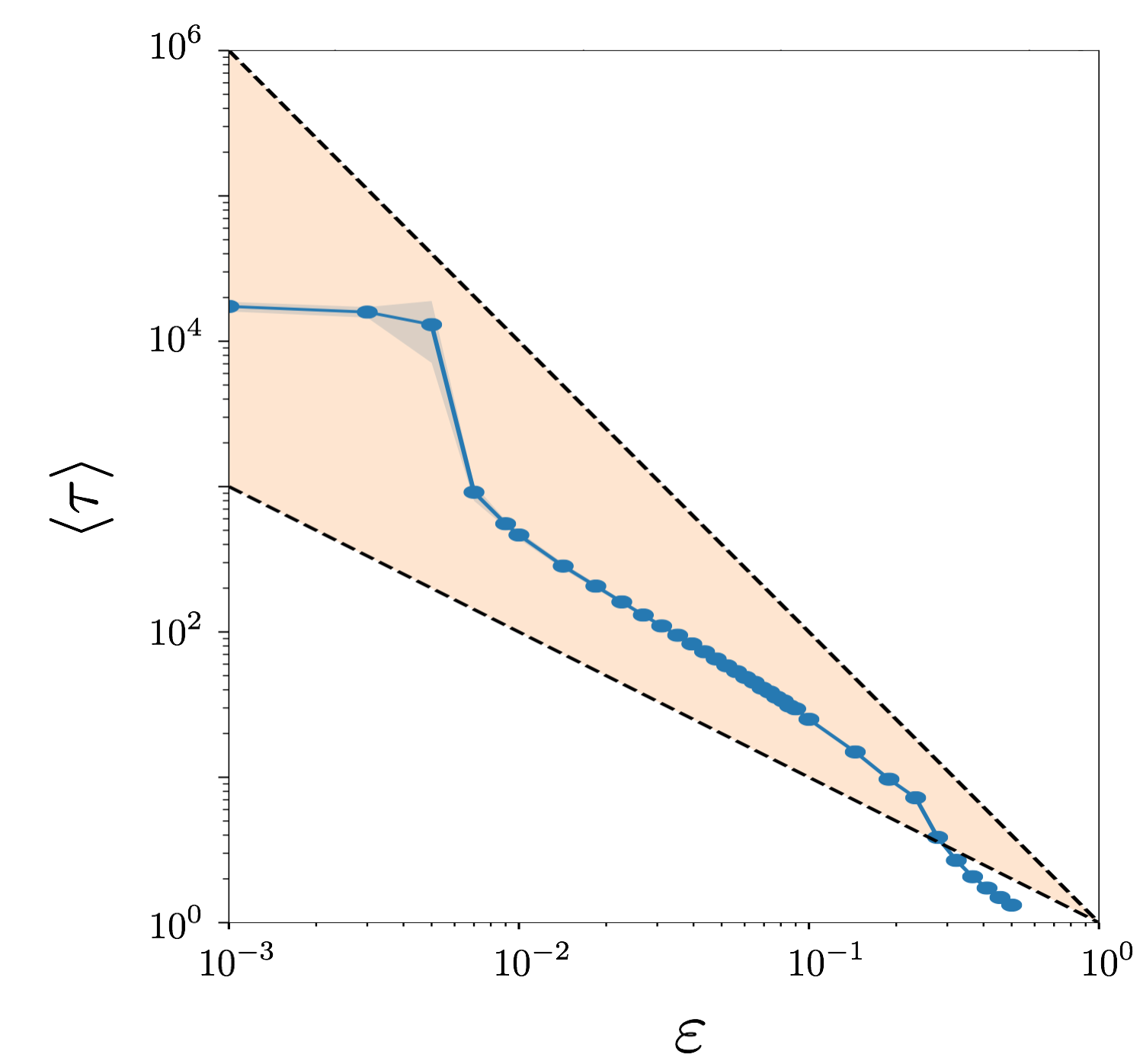}
\caption{The lifetime $\langle \tau \rangle$ as function of inter-coupling strength $\varepsilon$. The (dots) correspond to average over fifty initial conditions chosen nearby $M$ and (blue shaded) area is the standard deviation. The (dashed) lines represent the estimated scaling $\varepsilon^{-1}$ and $\varepsilon^{-2}$. Panel a) displays the case of coupling function does not depend on $\delta$ and panel b) include it. The parameters are $\beta = 200$, $\sigma = 1.2$ and $\eta = 0.25$. }
\label{fig:breaking_collapsing_scaling}
\end{figure}

The lifetime depends on the initial conditions and parameters of the network.  We select initial conditions starting in a neighborhood of $M$~\footnote{We select initial conditions satisfying the restrictions for system be in $M$ and we add a small random
uniform vector with each coordinate drawn from the interval $(0,0.01)$. }. And for each value of $\varepsilon$, we compute $\tau$ as shown in Figure~\ref{fig:breaking_collapsing_scaling}. The behavior of $\langle \tau \rangle$ as a function of $\varepsilon$ does not entirely seems to plateau after two decades. In interesting question is whether this plateau really exist or is a numerical artefact. As we only provide an lower bound for the survival time, this issue remains open.  To obtain description of such behavior we need to take a closer look at the dynamics near the boundary of the invariant manifold. As the dynamics hits the boundary it becomes high-dimensional which turns description of this situation challenging.

\vspace{0.5cm}
\noindent
{\bf Acknowledgements} Jaap Eldering was supported by FAPESP grant 2015/25947-6. Ed Roque dos Santos was supported by FAPESP grant 2018/10349-4, CAPES. TP is a Newton Advanced Fellow of the Royal Society NAF$\backslash$R1$\backslash$180236 and is also supported 
by Serrapilheira Institute (Grant No. Serra-1709-16124), FAPESP (grant 2013/07375-0)

\section*{References}

\bibliographystyle{amsplain}
\bibliography{chimeras}

\providecommand{\bysame}{\leavevmode\hbox to3em{\hrulefill}\thinspace}
\providecommand{\MR}{\relax\ifhmode\unskip\space\fi MR }
\providecommand{\MRhref}[2]{%
  \href{http://www.ams.org/mathscinet-getitem?mr=#1}{#2}
}
\providecommand{\href}[2]{#2}
\begin{thebibliography}{10}

\bibitem{abrams2008}
D~M Abrams, R~Mirollo, S~H Strogatz, and D~A Wiley, \emph{Solvable model for
  chimera states of coupled oscillators}, Phys. Rev. Lett. \textbf{101} (2008),
  084103.

\bibitem{jitcxde}
Gerrit Ansmann, \emph{Efficiently and easily integrating differential equations
  with {JiTCODE}, {JiTCDDE}, and {JiTCSDE}}, Chaos \textbf{28} (2018), 043116.

\bibitem{Ashwin2015}
P~Ashwin and O~Burylko, \emph{Weak chimeras in minimal networks of coupled
  phase oscillators}, Chaos \textbf{25} (2015), 013106.

\bibitem{barabasi1999}
A-L Barab{\'a}si and R~Albert, \emph{Emergence of scaling in random networks},
  Science \textbf{286} (1999), 509--512.

\bibitem{Bick2016}
C~Bick and P~Ashwin, \emph{Chaotic weak chimeras and their persistence in
  coupled populations of phase oscillators}, Nonlinearity \textbf{29} (2016),
  1468.

\bibitem{Chicone2006}
C~Chicone, \emph{Ordinary differential equations with applications}, 2nd ed.,
  Texts in applied mathematics, v.34, Springer, New York, 2006.

\bibitem{eldering2013}
J~Eldering, \emph{Normally hyperbolic invariant manifolds --- the noncompact
  case}, Atlantis Series in Dynamical Systems, vol.~2, Atlantis Press, Paris,
  September 2013.

\bibitem{fenichel1971}
N~Fenichel, \emph{Persistence and smoothness of invariant manifolds for flows},
  Indiana Univ. Math. J. \textbf{21} (1971/1972), 193--226.

\bibitem{Haugland2015}
S~W Haugland, L~Schmidt, and K~Krischer, \emph{Self-organized alternating
  chimera states in oscillatory media}, Scientific Reports \textbf{5} (2015),
  no.~9883.

\bibitem{HPS1977}
M~W Hirsch, C~C Pugh, and M~Shub, \emph{Invariant manifolds}, Lecture Notes in
  Mathematics, vol. 583, Springer-Verlag, Berlin, 1977.

\bibitem{Kementh_a_classification_2016}
F~P Kemeth, S~W Haugland, L~Schmidt, I~G Kevrekidis, and K~Krischer, \emph{A
  classification scheme for chimera states}, Chaos \textbf{26} (2016), 094815.

\bibitem{Ko_2008_2}
T-W Ko and G~B Ermentrout, \emph{Bistability between synchrony and incoherence
  in limit-cycle oscillators with coupling strength inhomogeneity}, Phys. Rev.
  E \textbf{78} (2008), 026210.

\bibitem{Ko_2008_1}
\bysame, \emph{Partially locked states in coupled oscillators due to
  inhomogeneous coupling}, Phys. Rev. E \textbf{78} (2008), 016203.

\bibitem{Kuramoto1984}
Y~Kuramoto, \emph{{Chemical Oscillations, Waves, and Turbulence}}, vol.~19,
  Springer Berlin Heidelberg, 1984.

\bibitem{kuramoto2002}
Y~Kuramoto and D~Battogtokh, \emph{Coexistence of coherence and incoherence in
  nonlocally coupled phase oscillators.}, Nonlinear Phenom. Complex Syst.
  (2002), no.~4, 380--385.

\bibitem{Larger2015}
L~Larger, B~Penkovsky, and Y~Maistrenko, \emph{Laser chimeras as a paradigm for
  multistable patterns in complex systems}, Nature Communications \textbf{6}
  (2015), no.~7752.

\bibitem{martens2013}
E~A Martens, S~Thutupalli, A~Fourri\`ere, and O~Hallatschek, \emph{Chimera
  states in mechanical oscillator networks}, Proc. Natl. Acad. Sci. USA
  \textbf{110} (2013), 10563--10567.

\bibitem{marvel2009}
S~A Marvel, R~E Mirollo, and S~H Strogatz, \emph{Identical phase oscillators
  with global sinusoidal coupling evolve by {M}\"obius group action}, Chaos
  \textbf{19} (2009), 043104, 11.

\bibitem{omel2013coherence}
O~E Omel'chenko, \emph{Coherence--incoherence patterns in a ring of non-locally
  coupled phase oscillators}, Nonlinearity \textbf{26} (2013), 2469.

\bibitem{Omel_chenko_2018}
\bysame, \emph{The mathematics behind chimera states}, Nonlinearity \textbf{31}
  (2018), R121--R164.

\bibitem{panaggio2015}
M~J Panaggio and D~M Abrams, \emph{Chimera states: coexistence of coherence and
  incoherence in networks of coupled oscillators}, Nonlinearity \textbf{28}
  (2015), R67--R87.

\bibitem{Sanders2007Averaging}
J~A Sanders, F~Verhulst, and J~Murdock, \emph{Averaging methods in nonlinear
  dynamical systems}, Applied Mathematical Sciences, Springer New York, 2007.

\bibitem{RalfPaper}
R~Toenjes, C~Fiore, and T~Pereira, \emph{Network induced coherence resonance},
  to appear (2020).

\bibitem{Totz_2018_spiral_wave_chimeras}
J~F Totz, J~Rode, M~R Tinsley, K~Showalter, and H~Engel, \emph{Spiral wave
  chimera states in large populations of coupled chemical oscillators}, Nature
  Physics \textbf{14} (2018), 282--285.

\bibitem{Vlasov_kuramoto_japanese_drums_2015}
V~Vlasov, A~Pikovsky, and E~E~N Macau, \emph{Star-type oscillatory networks
  with generic kuramoto-type coupling: A model for ``japanese drums
  synchrony''}, Chaos \textbf{25} (2015), 123120.

\bibitem{Vlasov2015}
V~Vlasov, Y~Zou, and T~Pereira, \emph{Explosive synchronization is
  discontinuous}, Phys. Rev. E \textbf{92} (2015), no.~1, 012904.

\bibitem{watanabe1994}
S~Watanabe and S~H Strogatz, \emph{Constants of motion for superconducting
  josephson arrays}, Phys. D \textbf{74} (1994), 197--253.

\bibitem{MATIAS}
M~Wolfrum and E~Omel'chenko, \emph{Chimera states are chaotic transients},
  Physical Review E \textbf{84} (2011), 015201.

\bibitem{wolfrum2011spectral}
M~Wolfrum, O~E Omel'chenko, S~Yanchuk, and Y~L Maistrenko, \emph{Spectral
  properties of chimera states}, Chaos \textbf{21} (2011), 013112.

\end{thebibliography}

\end{document}